\documentclass[a4paper,12pt]{article}
\usepackage{amsfonts}
\usepackage{amsmath}
\usepackage{amssymb}
\usepackage{graphicx}
\usepackage{epsf}
\usepackage{epsfig}

\usepackage{color}
\usepackage{afterpage}
\usepackage{verbatim}
\usepackage{amssymb}
\usepackage{graphicx}
\usepackage{amsmath}
\usepackage{amssymb}

\usepackage{graphicx}
\usepackage{epsfig}
\usepackage{amsmath}
\usepackage{color}
\usepackage{afterpage}
\usepackage{verbatim}

\usepackage[small]{caption}

\newtheorem{theorem}{Theorem}[section]
\newcommand{\R}{\mathbb{R}}

\newcommand{\VV}{\mathcal{V}}
\newcommand{\norm}[1]{\left\Vert#1\right\Vert}
\newcommand{\abs}[1]{\left\vert#1\right\vert}
\newcommand{\set}[1]{\left\{#1\right\}}
\newcommand{\para}[1]{\left(#1\right)}

\newcommand{\dive}{\textrm{div\,}}
\newcommand{\curl}{\textrm{curl\,}}
\newcommand{\D}{\mathrm{D}}
\newcommand{\B}{\mathrm{B}}

\begin{document}

\title{Optimization approach for the simultaneous 
reconstruction of the dielectric permittivity and magnetic
  permeability functions from limited observations}

\author{L. Beilina  \thanks{
Department of Mathematical Sciences, Chalmers University of Technology and
Gothenburg University, SE-42196 Gothenburg, Sweden, e-mail: \texttt{\
larisa@chalmers.se}}
 \and    Michel Cristofol \thanks{
Institut de Math\'{e}matiques de Marseille, Aix-Marseille Universit\'{e}
 13453 Marseille, France, e-mail: \texttt{\
michel.cristofol@univ-amu.fr}}
\and  Kati Niinim\"aki \thanks{ 
Department of Mathematical Sciences, KTH Royal Institute of Technology,
  SE-100 44 Stockholm, Sweden, e-mail: \texttt{niin@kth.se}}}

\maketitle

\graphicspath{   
{FIGURES/}
}

\begin{abstract}
We consider the inverse problem of the simultaneous reconstruction of the dielectric
permittivity and magnetic permeability functions of the Maxwell's system in
3D with limited boundary observations of the electric field. The
theoretical stability for the problem is provided by the Carleman estimates.
For the numerical computations the problem is formulated as an
optimization problem and hybrid finite element/difference method is
used to solve the parameter identification problem.
\end{abstract}

\section{Introduction}
  
This paper is focused on the numerical reconstruction of the
dielectric permittivity coefficient $\varepsilon(x)$ and the magnetic
permeability coefficient $\mu(x)$ for Maxwell's system basing our
observation on a single measurement data of the electric field
$E(x,t)$.
That means that we use boundary measurements of $E(x,t)$ which are
generated by a single direction of a plane wave.
For the numerical
reconstructions of the dielectric permittivity coefficient
$\varepsilon(x)$ and the magnetic permeability coefficient $\mu(x)$,
we consider a similar hybrid finite element/difference method (FE/FDM)
as was developed in \cite{BMaxwell}.

In the literature, the stability results for Maxwell's equations have
been proposed by Dirichlet-to-Neumann map or by Carleman estimates. For results using
Dirichlet-to-Neu-mann map with infinitely many boundary observations, see
\cite{BI,C,COS,KLS,OPS,SKU}. For results with a finite number of
observations several Carleman estimates have been derived in
\cite{BCS2012,EINT, Klibanov1,L05,LY1,LY2}. In order to solve the coefficient
inverse problem of Maxwell's equations numerically, we consider only
those theoretical results that involve finite number of observations.

Since real observations are generally corrupted by noise, it is
important to verify whether close observations lead to close
estimations of the coefficients. The results of this paper give such
conditions on the observations for estimating the parameters $\mu(x)$
and $\varepsilon(x)$ of the Maxwell's system (\ref{eq:maxwell}) from the
observations of $E(x,t)$ on the boundary of the domain. In
particular, we obtain a stability inequality of the form
$$\|(\mu_1,\varepsilon_1)-(\mu_2,\varepsilon_2)\|\le f\left ( \,
\|\hbox{Observation}(E_1)-\hbox{Observation}(E_2)\| \right )$$
$$
\mbox{for } f(s) \mbox{ such that } \displaystyle \lim_{s \rightarrow 0} f(s) =0,
$$ which links the distance between two sets of coefficients with the
distance between two sets of boundary observations of the electric
field $E(x,t)$. Such stability inequalities lead to the uniqueness of
the coefficients $(\mu(x),\varepsilon(x))$ given the observation
$"\hbox{Observation}(E(x,t))"$ on a small neighborhood of the boundary
of the domain of interest. They are also useful for the numerical
reconstruction of the coefficients using noise-free
observations~\cite{VdHQiu12}. Our main theoretical result concerns
stability inequality which gives an estimate of the norm of two
coefficients $\varepsilon$ and $\mu$ in terms of observation of only
the electric field $E(x,t)$ on the boundary of the domain. This
implies directly a uniqueness result. In the domain of the inverse
 problems, associated to reconstructing the Maxwell's coefficients from
 finite number of observations, except the reference \cite{Klibanov1}, and up to our knowledge, there exists no
 result involving only one component $E(x,t)$ (or the component
 $H(x,t)$).\\

In this work, the minimization problem of reconstructing functions
$\varepsilon(x)$ and $\mu(x)$ is reformulated as the problem of
finding a stationary point of a Lagrangian involving a forward
equation (the state equation), a backward equation (the adjoint
equation) and two equations expressing that the gradients with respect
to the coefficients $\varepsilon(x)$ and $\mu(x)$ vanish.  Moreover,
in our work the forward and adjoint problems are given by
time-dependent Maxwell's equations for the electric field. This means
that we have observations of the electric field in space and in time
 which provides us better reconstruction of both coefficients.
In order to get the computed values of $\varepsilon(x) $ and $\mu(x)$,
we arrange an iterative process by solving in each step the forward
and backward equations and updating the coefficients $\varepsilon(x) $
and $\mu(x)$ at every step of our iterations.

Recall, that in our optimization procedure the forward and adjoint
problems are given by the time-dependent Maxwell's system for the
electric field.  For the numerical solution of the Maxwell equations,
different formulations are available. One can consider, for example,
the edge elements of N\'ed\'elec \cite{Nedelec86}, the node-based
first-order formulation of Lee and Madsen \cite{LM}, the node-based
curl-curl formulation with divergence condition of Paulsen and Lynch
\cite{PL} or the interior-penalty discontinuous Galerkin FEM
\cite{GSS}.  
 In this work, for the discretization of the Maxwell's equations, we
 use stabilized domain decomposition method of \cite{BMaxwell} with
 divergence condition of Paulsen and Lynch \cite{PL} which removes
 spurious solutions when local mesh refinement is applied and material
 discontinuities are not too big \cite{BG}. Numerical tests of
 \cite{BMaxwell} in the case of coefficient inverse problems (CIPs), similar to our of recent
 experimental work \cite{BTKB}, show that these spurious solutions will not appear.

It is well known that edge elements are the most satisfactory from a
theoretical point of view ~\cite{Monk} since they automatically
satisfy the divergence free condition.  However, in the case of
time-dependent computations they are less attractive. First, they are
time consuming because the solution of a resulted linear system is
required at every time iteration. Second, in \cite{delta} was shown
that in the case of triangular or tetrahedral edge elements, the
entries of the diagonal matrix resulting from mass-lumping are not
necessarily strictly positive, and thus, explicit time stepping cannot
be used in general. The method of \cite{BMaxwell} where nodal elements
were used leads to a fully explicit and efficient scheme when mass-lumping is
applied \cite{joly,delta}. This method is efficiently
implemented in the software package WavES \cite{waves} in C++ using
PETSc \cite{petsc} and MPI (message passing interface) and 
 is convenient   for our simulations.

We note, that we reconstruct  functions $\varepsilon(x)$ and $\mu(x)$
simultaneously. Up to our knowledge, there are few studies (e.g. \cite{FJZ10}) that
consider the simultaneous recovery of $\varepsilon(x)$ and $\mu(x)$ within the
spatio-temporal Maxwell's equations. In the papers \cite{BTKB, BMaxwell}, which use the same
optimization approach, the coefficient $\mu$ is assumed to have
a known and a constant value, i.e. $\mu=1$, which means that the medium is
non-magnetic. However, it is well known that in many situations,
we need to deal with magnetic materials (e.g. metamaterials \cite{SSMS}, low-loss materials 
 \cite{PCP}).

Potential applications of our algorithm are in reconstructing the electromagnetic
parameters in nanocomposites or artificial materials \cite{PCP, SSD, SSMS}, imaging
of defects and their sizes in a non-destructive testing of materials
 and in photonic
crystals \cite{CKKJ}, measurement of the moisture content
\cite{CCSR01} and drying processes \cite{MDBSR01}, for example.

Our numerical simulations show that we are able to accurately
reconstruct simultaneously contrasts for both functions
$\varepsilon(x)$ and $\mu(x)$ as well as their locations. In our
future work, similarly with \cite{BTKB, BMaxwell, BJ}, we are planning
also to reconstruct shapes of the inclusions using a posteriori error
estimates in the Tikhonov functional or in the reconstructed
coefficients and based on them an adaptive finite element method.

The paper is organized as follows. In Section 2 we state the
theoretical results of recovering $\varepsilon(x)$ and $\mu(x)$ from
the limited boundary measurements of $E(x,t)$. Section 3 is devoted to
presenting the numerical method used in this article. In Section 4, we
give detailed information about our discrete numerical method and we
outline the algorithm for the solution of our
inverse problem. The
numerical results are presented in Section 5. Discussions and
conclusions are given in Section 6.

\section{Theoretical results}

Let us consider a bounded domain $\Omega \subset \mathbb{R}^3$  with a smooth
boundary $\partial \Omega$, and define by $\Omega_T := \Omega \times
(0,T)$, $\partial \Omega_T := \partial \Omega \times (0,T)$,  where
$T$ is a strictly positive constant.
The electromagnetic equations in an inhomogeneous isotropic case in the bounded domain  $\Omega \subset  \mathbb{R}^d,~ d=2,3$  with boundary $\partial \Omega$,  are described by the first order system of partial differential equations
\begin{equation}
  \label{eq:maxwell}
  \left \{ \begin{array}{llllll} 
\partial_{t} D - \nabla \times H(x,t) = 0 &&\mbox{ in } \Omega_T\\
 \partial_{t} B + \nabla \times E(x,t) = 0, &&
 \mbox{ in } \Omega_T,\\
 D(x,t)= \varepsilon E(x,t), \quad  B(x,t)= \mu H(x,t),&&\\
E(x,0) = E_0(x), \quad  H(x,0) = H_0(x), &&\\
\nabla \cdot D(x,t) = 0,\quad \nabla \cdot B(x,t)  =0, && \mbox{ in } \Omega_T,\\
         D(x,t)\times\nu=0,\quad B(x,t)\cdot\nu=0, && \textrm{on }\,\partial \Omega_T, \end{array} \right .
\end{equation}
where  $E(x,t), H(x,t), D(x,t), B(x,t)$ are three-dimensional
vector-valued functions of the time $t$ and the space variable
$x=(x_1,x_2,x_3)$, and correspond to  the electric and magnetic
fields  and the electric and magnetic inductions, respectively. 
The  dielectric permittivity, $\varepsilon(x)>0$  and the magnetic permeability, $\mu(x)>0$, depend on   $x \in \Omega$,  $\nu=\nu(x)$ denotes the unit outward normal vector to $\partial \Omega$.

Our goal is to reconstruct the coefficients $\varepsilon(x)$ and $\mu(x)$ in the
system (\ref{eq:maxwell}) with appropriate initial 
conditions $E_{0}$ and $H_{0}$ on the electric and magnetic inductions, using only a finite number of observations of the electric
field $E(x,t)$ on the boundary $\partial \Omega$ of the domain $\Omega$.

We base our theoretical approach on the work  \cite{BCS2012}. In this
work the authors obtained a stability inequality for the dielectric
permittivity $\varepsilon(x)$ and the magnetic permeability $\mu(x)$ involving
boundary observations of both magnetic induction $B(x,t)$ and electric
induction $D(x,t)$ (see theorem 1 of \cite{BCS2012}). 
With the same assumptions, as those in \cite{BCS2012}, we can reformulate this stability
result using only the observations of the electric field $E(x,t)$ (or
correspondingly, using
only the observations of the magnetic field $H(x,t)$):

Assume that the functions $\mu(x)$ and $\varepsilon(x)$ in
$\mathcal{C}^2(\overline{\Omega}), x\in\overline{\Omega}$ obey
\begin{equation}\label{1.2}
\mu(x) \geq \mu_0, \quad   \varepsilon(x) \geq \varepsilon_0,
\end{equation}
for some $\varepsilon_0>0$ and $\mu_0>0$. Next, for simplicity,  we  introduce some
similar notations to the ones introduced in \cite{BCS2012}. 

Let us pick $x_0\in\R^3\backslash\overline{\Omega}$, set $c(x)=(\mu(x) \varepsilon(x))^{-1}$ for $x \in \overline{\Omega}$, $c_0=(\mu_0 \varepsilon_0)^{-1}$ and assume that the following condition holds for some $\rho \in (0,c_0)$
\begin{equation}\label{1.6}
\frac{3}{2}\abs{\nabla\log c(x)}\abs{x-x_0}\leq 1-\frac{\rho}{c_0}, \quad x\in\overline{\Omega}.
\end{equation}
This  technical condition is claimed by the weight function $\psi_0
=|x-x_0|^2$ that is used to create the Carleman estimate established
to prove the theorem \ref{T.1} in \cite{BCS2012}. In other terms,  (\ref{1.6}) arises
from the classical pseudo-convexity condition. Another standard
hypothesis is that the coefficients $\varepsilon$ and $\mu$ are known
in a neighborhood of the boundary of $\Omega$ 

Next, we define $\omega=\Omega\cap \mathcal{O}$ where $\mathcal{O}$ is some
neighborhood of $\partial \Omega$ in $\R^3$. Further let
$M_0>0$ and two given functions $\mu^\sharp,\lambda^\sharp$ belong to
$\mathcal{C}^2(\omega)$. Now we can define the admissible set of unknown coefficients $\mu$ and $\varepsilon$ as
\begin{equation}\label{1.7}
\Lambda_\omega(M_0)=\set{(\mu,\varepsilon)\ \textrm{verifying}\,(\ref{1.2}) ;
\norm{(\mu,\varepsilon)}_{\mathcal{C}^2(\overline{\Omega})}\leq M_0\ \textrm{and}
\,(\mu,\varepsilon)=(\mu^\sharp,\varepsilon^\sharp)\ \textrm{in}\ \omega}.\end{equation}
We set 
$$\VV= H_{\tau,0}(\curl,\dive 0;\Omega) \times {\rm H}_{n,0}(\curl,\dive 0;\Omega), $$
where
$$H_{*,0}(\curl,\dive 0;\Omega)=\{ u \in H^1(\Omega)^3,\ \dive u = 0\ {\rm and}\ \gamma_* u = 0 \},\ *=\tau, n.$$

Further, for the identification of $(\mu,\varepsilon)$,  imposing  (as
will appear in the sequel) that $(B(x,t),D(x,t))$ are observed twice,
we consider two sets of initial data $(D_0^k,B_0^k)$, $k=1,2$ such that,
\begin{equation*}
\D_0^k(x)=\para{d_1^k(x),d_2^k(x),d_3^k(x)}^\top,\quad\B_0^k(x)=\para{b_1^k(x),b_2^k(x),b_3^k(x)}^\top,
\end{equation*}
and define the $12\times 6$ matrix
\begin{equation}\label{1.9}
\mathcal{K}(x)=\left(
              \begin{array}{cccccc}
                e_1\times\B_0^1 &  e_2\times\B_0^1 &  e_3\times\B_0^1 & 0 & 0 & 0 \\
                0 & 0 & 0 &  e_1\times\D_0^1 &  e_2\times\D_0^1 &  e_3\times\D_0^1 \\
                 e_1\times\B_0^2 &  e_2\times\B_0^2 &  e_3\times\B_0^2 & 0 & 0 & 0 \\
                0 & 0 & 0 &  e_1\times\D_0^2 &  e_2\times\D_0^2 &  e_3\times\D_0^2 \\
              \end{array}
            \right),
\end{equation}
where $ x\in\Omega$. Then we can write that $(\B_i^k(x,t),\D_i^k(x,t))$ are the solution to
(\ref{eq:maxwell}) with the initial data $(\B_0^k,\D_0^k)$, $k=1,2$, where
$(\mu,\varepsilon)$ are substituted with $(\mu_i,\varepsilon_i)$, $i=1,2$.\\
Finally, we note that
$\mathcal{H}( \Omega_{T})=H^3(-T,T;L^2(\Gamma))\cap
H^2(-T,T;H^1(\Gamma))$ is a  Hilbert space 
equipped with the norm
$$
\norm{u}^2_{\mathcal{H}( \Omega_{T})}=\norm{u}^2_{H^3(-T,T;L^2(\Gamma))}+\norm{u}^2_{H^2(-T,T;H^1(\Gamma))},\quad u\in \mathcal{H}( \Omega_{T}).
$$
The extension of the time interval to (-T,T) corresponds to a technical point in the proof of the stability inequality  (see lemma 3 in \cite{BCS2012}).\\
Now we recall the main theoretical result of the paper \cite{BCS2012}:
Under some hypothesis on $T$ and choosing initial
conditions  $(\B_0^k,\D_0^k)$, $k = 1,2$ verifying some additional
assumptions, then there are two constants $C>0$ and $\kappa\in (0,1)$, depending on $\Omega$, $\omega$, $T$, $M$ and $M_0$, such that we have:
\begin{eqnarray*}
\norm{\mu_1-\mu_2}_{H^2(\Omega)}&+&\norm{\varepsilon_1-\varepsilon_2}_{H^2(\Omega)}\\
&&   \leq C\para{\sum_{k=1}^2\para{\norm{\para{\B_1^k-\B_2^k}_\tau}_{\mathcal{H}( \Omega_{T})}+ \norm{\para{\D_1^k-\D_2^k}_\nu}_{\mathcal{H}( \Omega_{T})}}}^\kappa.
\end{eqnarray*}

Under the same  assumptions and considering the definition of the electric induction
$\D_i$ and the magnetic induction  $\B_i$  for $i=1,2$ we can write in the neighborhood $\omega$ of the boundary $\partial \Omega$ the following relations: 
$$
\partial_{t} B_{i}^{k}= -\curl \varepsilon_i^{-1} \D_{i}^{k} , \mbox{ in } \omega \times (-T,T), \mbox{ for } i=1,2 \mbox{ and } k=1,2.
$$
 Since $\varepsilon_i = \varepsilon^\sharp$, for  $i=1,2$ in $\omega$, we can write 
$$
\partial_{t}(\B^k_1-\B^k_2) = -\curl \frac{1}{\varepsilon^\sharp} (\D_1^k-\D_2^k) \mbox{ in } \omega \times (-T,T), \mbox{ for } k=1,2.
$$
It is straightforward to verify that
$$
\partial_{t}(\B^k_1-\B^k_2) = - \frac{1}{\varepsilon^\sharp} \curl  (\D_1^k-\D_2^k) - \nabla \frac{1}{\varepsilon^\sharp} \times (\D_1^k-\D_2^k).
$$
We define by
$$
N(\partial \Omega_{T})=H^3(-T,T;L^2(\partial \Omega))\cap
H^2(-T,T;H^1(\partial \Omega))\cap H^1(-T,T;H^2(\partial \Omega)),
$$
the Hilbert space 
equipped with the norm
$$
\norm{u}^2_{N(\partial
  \Omega_{T})}=\norm{u}^2_{H^3(-T,T;L^2(\partial
  \Omega))}+\norm{u}^2_{H^2(-T,T;H^1(\partial \Omega))}+\norm{u}^2_{H^1(-T,T;H^2(\partial \Omega))},\; 
$$
$u\in N(\partial \Omega_{T}).$Then from \eqref{1.7} and since $B_{1}^{k}(0,\cdot) =
B_{2}^{k}(0,\cdot)$ and $D_{1}^{k}(0,\cdot) = D_{2}^{k}(0,\cdot)$, for
$k=1,2$, we get
$$
\norm {\B^k_1-\B^k_2}_{\mathcal{H}( \Omega_{T})} \leq M_0 \norm{\D^k_1-\D^k_2}_{N(\partial \Omega_{T})}. 
$$
Thus we can deduce our theorem in the following form

\begin{theorem} \label{T.1} Let $T>c_0^{-1 \slash 2} \max_{x\in\overline{\Omega}}\abs{x-x_0}$ and pick $(\B_0^k,\D_0^k)\in (H^2(\Omega)^3 \times H^2(\Omega)^3) \cap \VV$, $k=1,2$, in such a way that there exists a $6 \times 6$ minor $m(x)$ of
the matrix $\mathcal{K}(x)$ defined in \eqref{1.9}, obeying:
\begin{equation*}
m(x) \neq 0,\ x \in \overline{\Omega\backslash\omega}.
\end{equation*}
Further, choose $(\mu_i,\varepsilon_i)\in\Lambda_\omega(M_0)$, $i=1,2$, such that
\begin{equation*}
\norm{\para{\B_i^k,\D_i^k}}_{\mathcal{C}^3(-T,T;W^{2,\infty}(\Omega))}\leq M,\ k=1,2,
\end{equation*}
for some $M>0$. Then there are two constants $C>0$ and $\kappa\in (0,1)$, depending on $\Omega$, $\omega$, $T$, $M$ and $M_0$, such that we have:
$$
\norm{\mu_1-\mu_2}_{H^2(\Omega)}+\norm{\varepsilon_1-\varepsilon_2}_{H^2(\Omega)}\leq
C\left
  ({\sum_{k=1}^2{\norm{{D_1^k-D_2^k}_\nu}_{N(\partial \Omega_{T})}}}\right
)^\kappa.
$$
\end{theorem}

\section{Statement of the forward and inverse problems}
\subsection{The mathematical model}
Below for any vector function $ u \in \mathbb{R}^3$ our notations $u
\in L^2(\Omega)$ or $u \in H^k(\Omega), k=1,2$, mean that every
component of the vector function $u$ belongs to this space.

\begin{figure}[tbp]
 \begin{center}
 \begin{tabular}{c}
 {\includegraphics[width=12.0cm, clip = true, trim = 0.0cm 0.0cm 0.0cm 0.0cm]{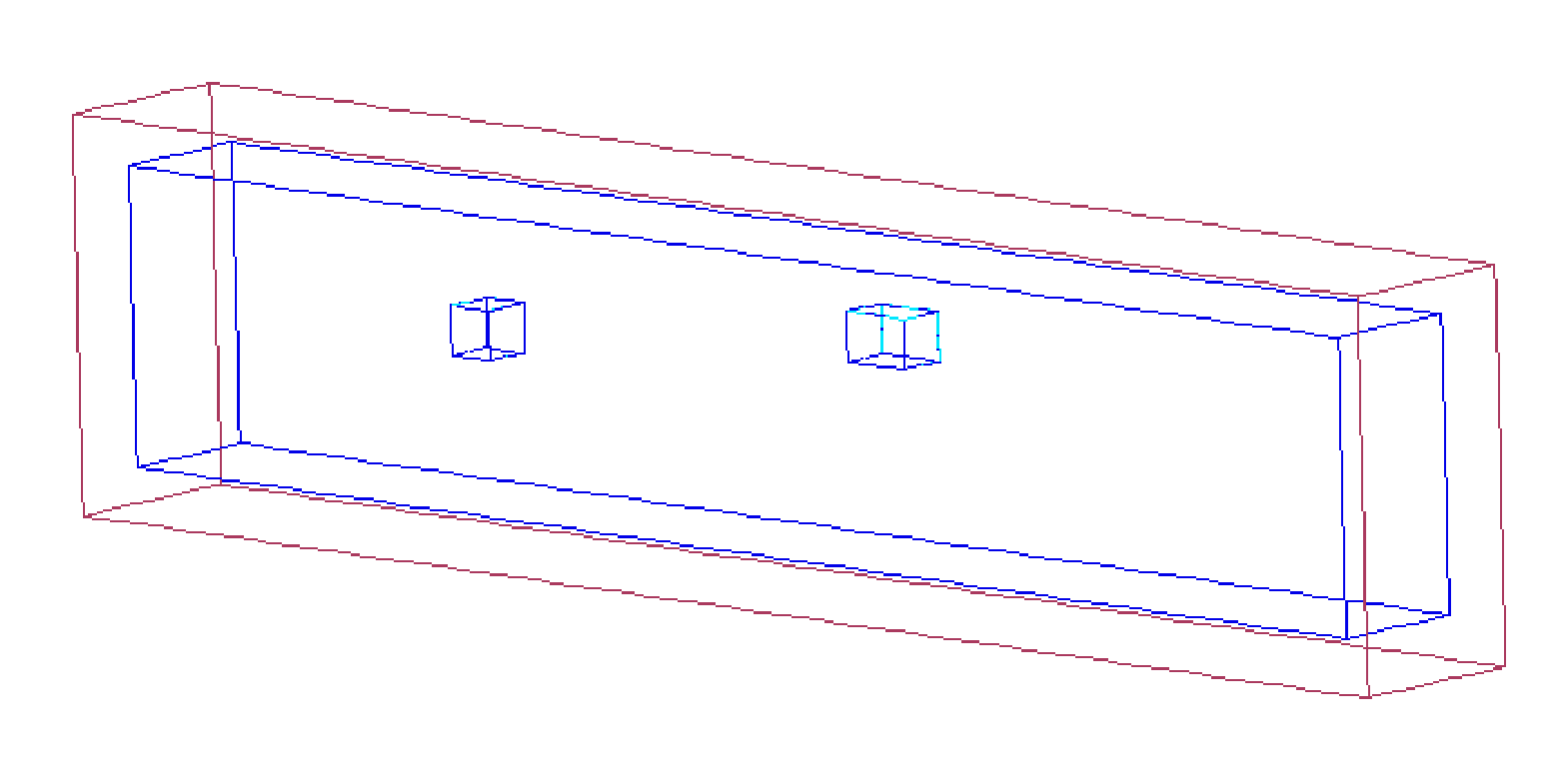}} \\
 a) $\Omega = \Omega_{\rm FEM} \cup \Omega_{\rm FDM}$ \\
  {\includegraphics[width=12.0cm, clip = true, trim = 0.0cm 0.0cm 0.0cm 0.0cm]{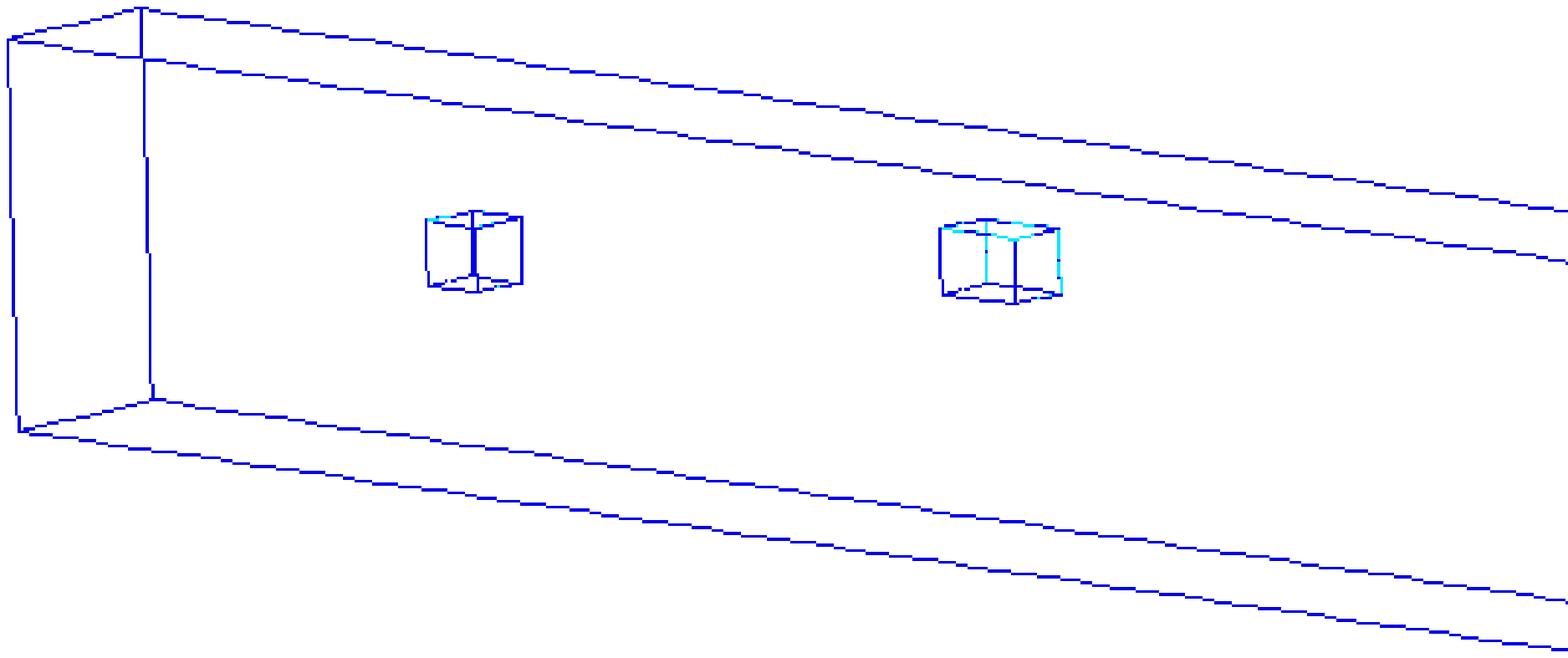}} \\
  b)  $\Omega_{\rm FEM}$ 
 \end{tabular}
 \end{center}
 \caption{{\protect\small \emph{Illustration of the domain $\Omega$. a)
       Combination of domains  $\Omega_{\rm FDM}$ and $\Omega_{\rm FEM}$. b)
       Only the domain $\Omega_{\rm FEM}$, with two inclusions. }}}
 \label{fig:fig1}
 \end{figure}

Next, we decompose
 $\Omega$ into two subregions, $\Omega_{\rm FEM}$ and $\Omega_{\rm FDM}$ such
 that $\Omega = \Omega_{\rm FEM} \cup \Omega_{\rm FDM}$,  $\Omega_{\rm FEM}
 \cap \Omega_{\rm FDM} = \emptyset$ and $\partial \Omega_{\rm FEM}
 \subset \partial \Omega_{\rm FDM}$, for an illustration of the domain
 decomposition, see figure \ref{fig:fig1}. 
 In $\Omega_{\rm FEM}$ we use finite elements and in  $\Omega_{\rm FDM}$ we will use finite difference method with first order
 absorbing boundary conditions \cite{EM}. The boundary $\partial \Omega$ is  such that $\partial \Omega
=\partial _{1} \Omega \cup \partial _{2} \Omega \cup \partial _{3} \Omega$  where
$\partial _{1} \Omega$ and $\partial _{2} \Omega$ are, respectively, front and
back sides of the domain $\Omega$, and $\partial _{3} \Omega$ is the union
of left, right, top and bottom sides of this domain.

By eliminating $B$ and $D$ from (\ref{eq:maxwell}) we obtain
the model problem for the electric field $E$ with the
perfectly conducting boundary conditions
 at the boundary
$\partial \Omega$ is as follows:
\begin{eqnarray}
\varepsilon \frac{\partial^2 E}{\partial t^2} + \nabla \times ( \mu^{-1} \nabla \times E)    &=& 0,~ \mbox{in}~~ \Omega_T,    \label{model1_1} \\
\nabla \cdot (\varepsilon E) &=& 0,~ \mbox{in}~~ \Omega_T,  \label{model1_2}  \\
  E(x,0) = f_0(x), ~~~E_t(x,0) &=& f_1(x)~ \mbox{in}~~ \Omega,      \label{model1_3}  \\
E \times n &=& 0 ~ \mbox{on}~~ \partial \Omega_T.  \label{model1_4} 
\end{eqnarray}
Here we assume that
\begin{equation*}
f_{0}\in H^{1}(\Omega), f_{1}\in L_{2}(\Omega). 
\end{equation*}
We note that a similar equation can be derived also for $H$.
 For numerical solution of (\ref{model1_1})-(\ref{model1_4}) in $\Omega_{\rm FDM}$ we  use the
finite difference method on a structured mesh with constant
coefficients $\varepsilon= \varepsilon^\sharp\ = 1$ and $\mu =
\mu^{\sharp} =1$. In $\Omega_{\rm FEM}$, we use finite elements on a sequence of
unstructured meshes $K_h = \{K\}$, with elements $K$ consisting of
triangles in $\mathbb{R}^2$ and tetrahedra in $\mathbb{R}^3$
satisfying maximal angle condition \cite{Brenner}.

In this work, for the discretization of the Maxwell's equations we use
stabilized domain decomposition method of \cite{BMaxwell} and consider
Maxwell's system in a convex geometry without reentrant corners where
material discontinuities are not too big \cite{BG}.  Since in our
numerical simulations the relative permeability $\varepsilon$ and
relative permittivity $\mu$ does not vary much, such assumptions about
the coefficients are natural. 
Moreover, recent experimental works \cite{ BTKB,BTKF,NBKF} show that
material discontinuities in real applications of reconstruction of
dielectrics (or refractive indices) of unknown targets placed in the
air or ground using microwave imaging technology are not big.  For
application where the coefficients are assumed to be smooth see \cite{WSKVH09} and \cite{ZTRL11}. \cite{WSKVH09} considers breast imaging and \cite{ZTRL11} soil moisture imaging. In our work we treat discontinuities of
coefficients $\epsilon$ and $\mu$ in a similar way as in formula
(3.25) of
section 3.5 in  \cite{CWZ14}. See also \cite{BMaxwell2} for the case
when only $\epsilon$ is unknown. Further, in our computations all materials with values
of $\varepsilon > 10$ are treated as metals and we call $\varepsilon$ as
``apparent''  or ``effective''  dielectric constant'', see \cite{BTKB, BTKF, NBKF} for more
information and explanation.

To stabilize the finite element solution using  standard piecewise
continuous   functions, we enforce
the divergence condition (\ref{model1_2}) and add a Coulomb-type gauge
condition \cite{Ass,div_cor} to (\ref{model1_1})-(\ref{model1_4}) with
$ 0 \leq s \leq 1$
 \begin{eqnarray}
 \varepsilon \frac{\partial^2 E}{\partial t^2} + \nabla \times ( \mu^{-1} \nabla \times E)  - s\nabla  ( \nabla \cdot(\varepsilon E))  &=& 0,~ \mbox{in}~~ \Omega_T,    \label{model2_1} \\
   E(x,0) = f_0(x), ~~~E_t(x,0) &=& f_1(x)~ \mbox{in}~~ \Omega,      \label{model2_3}  \\
 E \times n &=& 0
~~ \partial \Omega_T,  \label{model2_4} \\
 \mu(x)=\varepsilon \left( x\right) &=&1\text{ in }\Omega _{\rm FDM}. \label{extra2}
 \end{eqnarray}

Let $S_T := \partial_1 \Omega \times (0,T)$ where
 $\partial_1 \Omega$ is the backscattering side of the domain $\Omega$ with
 the time domain observations, and define by $S_{1,1} := \partial_1
 \Omega \times (0,t_1]$, $S_{1,2} := \partial_1 \Omega \times
 (t_1,T)$,  $S_2 := \partial_2 \Omega \times (0, T)$,  $S_3 :=
 \partial_3 \Omega \times (0, T)$.
Our forward problem used in computations, thus writes
\begin{equation}\label{E_gauge}
\begin{split}
\varepsilon \frac{\partial^2 E}{\partial t^2} + \nabla \times ( \mu^{-1} \nabla \times E)  - s\nabla  ( \nabla \cdot(\varepsilon E))  &= 0,~ \mbox{in}~~ \Omega_T, \\
  E(x,0) = 0, ~~~E_t(x,0) &= 0~ \mbox{in}~~ \Omega,     \\
\partial _{n}E& = p\left( t\right) ,~\mbox{on}~ S_{1,1},
\\
\partial _{n} E& =-\partial _{t} E,~\mbox{on}~ S_{1,2},
\\
\partial _{n} E& =-\partial _{t} E,~\mbox{on}~ S_2, \\
\partial _{n} E& =0,~\mbox{on}~ S_3.\\
\end{split}
\end{equation}

We use the Neumann boundary conditions at the left and right hand
sides of a   domain $\Omega$  (recall, that $\Omega = \Omega_{\rm FEM} \cup \Omega_{\rm FDM}$  and in $\Omega_{\rm FDM}$  our coefficients $\varepsilon(x)=\mu(x)=1$)
 and first order absorbing boundary conditions
\cite{EM}  at the rest of the boundaries.
Application of Neumann boundary conditions allows us to assume infinite
structure in  lateral directions and thus, allows us to consider the
CIPs  of the reconstruction of unknown parameters $ \varepsilon$ and $\mu$  in a waveguide.
 Further, we assume
homogeneous initial conditions.

It was demonstrated numerically in \cite{BMaxwell} that the solution
of the problem (\ref{E_gauge}) approximates well the solution of the
original Maxwell's equations for $\mu=1, s=1.$ The energy estimate of
Theorem 4.1 of \cite{BMaxwell} implies an uniqueness result for the forward 
problem (14)  with $ \varepsilon \geq 1, \mu=1$.  Similar
theorem can be proven for the case when both coefficients are unknown,
and it can be done in a future research.

 We assume that our
 coefficients $\varepsilon \left(x\right), \mu(x) $ of equation (\ref{E_gauge})
 are such that
\begin{equation} \label{2.3}
\begin{split}
\varepsilon \left(x\right) &\in \left[ 1,d_1\right],~~ d_1 = const.>1,~ \varepsilon(x) =1
\text{ for }x\in  \Omega _{\rm FDM}, \\
\mu(x) &\in \left[ 1,d_2\right],~~ d_2 = const.>1,~ \mu(x) =1   \text{ for }x\in  \Omega _{\rm FDM}, \\
\varepsilon \left(x\right), \mu(x) &\in C^{2}\left( \mathbb{R}^{3}\right) . 
\end{split}
\end{equation}
The values of constants $d_1, d_2$ in
(\ref{2.3}) are chosen from real life experiments  similarly with  \cite{BTKB, BTKF,SSMS,NBKF} and we assume that
we know them a priori.

We consider the following

\textbf{Inverse Problem  (IP)} \emph{Suppose that the coefficients
}$\varepsilon\left(x\right)$ and $\mu(x)$ \emph{\ satisfies (\ref{2.3})
  such that  numbers }$d_1, d_2 >1$\emph{\ are given. Assume that the
  functions }$ \varepsilon\left(x\right), \mu(x) $\emph{\ are unknown in the
  domain }$\Omega \backslash  \Omega_{\rm FDM}$\emph{. Determine the functions }$ \varepsilon\left(x\right), \mu(x) $\emph{\ for }$x\in \Omega \backslash  \Omega_{\rm FDM},$ \emph{\ assuming that the
  following function }$\tilde E\left(x,t\right) $\emph{\ is known}
\begin{equation*}
 E \left(x,t\right) = \tilde E \left(x,t\right) ,\forall \left( x,t\right) \in
  S_T.  
\end{equation*}
A priori knowledge of an upper and lower bounds of functions
$\varepsilon\left(x\right)$ and $\mu(x)$ corresponds well with the
inverse problems concept about the availability of a priori
information for an ill-posed problem \cite{BKS, Engl, tikhonov}. In
applications, the assumption $\varepsilon\left(x\right) =\mu(x) =1$
for $x\in \Omega_{\rm FDM} $ means that the functions
$\varepsilon\left(x\right)$ and $\mu(x)$ have a known constant value
outside of the medium of interest $\Omega \backslash \Omega_{\rm FDM}
.$ The function $\tilde E\left(x,t\right) $ models time dependent
measurements of the electric wave field at the backscattering boundary
$\partial_1 \Omega$ of the domain of interest. In practice,
measurements are performed on a number of detectors, see \cite{BTKF,
  NBKF}.

\subsection{Optimization method}

We reformulate our inverse problem as an optimization problem, where
 we seek for two functions,  the permittivity $\varepsilon(x)$ and permeability $\mu(x)$, which result in a solution of
equations (\ref{E_gauge}) with best fit to time and space domain observations
$\tilde E$, measured at a finite number of observation points on
$\partial_1 \Omega$. 
Our goal is to minimize  the Tikhonov functional
\begin{equation}
F(E, \varepsilon, \mu) = \frac{1}{2} \int_{S_T}(E - \tilde{E})^2 z_{\delta }(t)
 d\sigma dt +
\frac{1}{2} \gamma_1  \int_{\Omega}( \varepsilon -  \varepsilon_0)^2~~ dx 
+ \frac{1}{2} \gamma_2  \int_{\Omega}( \mu -  \mu_0)^2~~ dx,
\label{functional}
\end{equation}
where $\tilde{E}$ is the observed electric field, $E$ satisfies the
equations (\ref{E_gauge}) and thus depends on $\varepsilon$ and $\mu$,
 $\varepsilon _{0}$ is the initial
guess for $\varepsilon $ and $\mu_{0}$ is the initial guess for $\mu$,
and $\gamma_i, i=1,2$ are the regularization parameters.  Here $z_{\delta
}(t)$ is a cut-off function, which is introduced to ensure that
the compatibility conditions at $\overline{\Omega}_{T}\cap \left\{
t=T\right\} $ for the adjoint problem (\ref{adjoint}) are satisfied, and $\delta
>0$ is a small number. We choose a function $z_{\delta }$
such that
\begin{equation*}
z_{\delta }\in C^{\infty }\left[ 0,T\right] ,\text{ }z_{\delta
}\left(t\right) =\left\{ 
\begin{array}{ll}
1 & \text{ for } t \in \left[ 0,T-\delta \right],  \\ 
0\ & \text{ for }t\in \left(T-\frac{\delta }{2},T\right],  \\ 
 0 < z_{\delta } < 1 & \text{ for }t\in \left(T-\delta ,T-\frac{\delta }{2}\right). 
\end{array}
\right.
\end{equation*}

Next, we introduce the following spaces of real valued vector functions
\begin{equation*}
\begin{split}
H_E^1 &:= \{ w \in H^1(\Omega_T):  w( \cdot , 0) = 0 \}, \\
H_{\lambda}^1 &:= \{ w \in  H^1(\Omega_T):  w( \cdot , T) = 0\},\\
U^{1} &=H_{E}^{1}(\Omega_T)\times H_{\lambda }^{1}(\Omega_T)\times C\left(\overline{\Omega}\right)\times C\left(\overline{\Omega}\right),\\
U^{0} &=L_{2}\left(\Omega_{T}\right) \times L_{2}\left(\Omega_{T}\right) \times
L_{2}\left(\Omega \right)\times
L_{2}\left(\Omega \right). 
\end{split}
\end{equation*}

To solve the minimization problem, we introduce the Lagrangian
\begin{equation*}
\begin{split}
L(u) = F(E, \varepsilon, \mu) 
&-  \int_{\Omega_T} \varepsilon  \frac{\partial
 \lambda }{\partial t} \frac{\partial E}{\partial t}  ~dxdt  
+   \int_{\Omega_T}( \mu^{-1}\nabla \times E)( \nabla \times \lambda)~dxdt  \\
&+ s  \int_{\Omega_T} ( \nabla \cdot (\varepsilon E))( \nabla \cdot \lambda)~dxdt \\
&- \int_{S_{1,1}} \lambda p(t) ~d \sigma dt    +
 \int_{S_{1,2}} \lambda \partial_t E ~d\sigma dt    
+\int_{S_2} \lambda \partial_t E  ~d\sigma dt, \\
\end{split}
\end{equation*}
where $u=(E,\lambda, \varepsilon, \mu) \in U^1$, and search for a stationary
point with respect to $u$ satisfying $ \forall \bar{u}= ( \bar{E}, \bar{\lambda},  \bar{\varepsilon}, \bar{\mu}) \in U^1$
\begin{equation}
 L'(u; \bar{u}) = 0 ,  \label{scalar_lagr}
\end{equation}
where $ L^\prime (u;\cdot )$ is the Jacobian of $L$ at $u$.

We assume that $\lambda \left(x,T\right) =\partial _{t}\lambda \left(x,T\right) =0$ and seek to impose such conditions on the function
$\lambda $ that $ L\left(E,\lambda ,\varepsilon, \mu \right)
:=L\left(u\right) =F\left(E,\varepsilon, \mu \right).$ Next, we use
the fact that $\lambda (x ,T) = \frac{\partial \lambda}{\partial t}
(x,T) =0$ and $E(x,0)= \frac{\partial E}{\partial t} (x ,0) = 0$, as
well as $\mu=\varepsilon=1$ on $\partial \Omega$, together with
boundary conditions $ \partial_n E = 0$ and $  \partial_n \lambda = 0$ on
$S_3$.  The equation (\ref{scalar_lagr}) expresses that
for all $\bar{u}$,
\begin{equation}\label{forward}
\begin{split}
0 = \frac{\partial L}{\partial \lambda}(u)(\bar{\lambda}) =
&- \int_{\Omega_T} \varepsilon \frac{\partial
\bar{\lambda}}{\partial t} \frac{\partial E}{\partial t}~ dxdt 
+  \int_{\Omega_T}  ( \mu^{-1} \nabla \times E)
(\nabla \times \bar{\lambda}) ~ dxdt  \\
&+ s  \int_{\Omega_T}(\nabla \cdot(\varepsilon E))( \nabla \cdot \bar{\lambda}) ~ dxdt \\
&- \int_{S_{1,1}} \bar{\lambda} p(t) ~d \sigma dt    + \int_{S_{1,2}}
\bar{\lambda} \partial_t E ~d\sigma dt   \\
& + \int_{S_2} \bar{\lambda} \partial_t E  ~d\sigma dt,~~\forall \bar{\lambda} \in H_{\lambda}^1(\Omega_T),\\
\end{split}
\end{equation}
\begin{equation} \label{control}
\begin{split}
0 = \frac{\partial L}{\partial E}(u)(\bar{E}) &=
\int_{S_T}(E - \tilde{E})~ \bar{E}~ z_{\delta}~ d \sigma dt \\
&-  \int_{\Omega_T} \varepsilon \frac{\partial
\lambda}{\partial t} \frac{\partial \bar{E}}{\partial t}~ dxdt
 + \int_{\Omega_T} (\mu^{-1}  \nabla \times \lambda)
(\nabla \times \bar{E})  ~ dxdt \\
&- \int_{\Omega} 
\frac{\partial{\lambda}}{\partial t}(x,0) \bar{E}(x,0) ~dx 
- \int_{S_{1,2 \cup S_2}} 
\frac{\partial{\lambda}}{\partial t} \bar{E} ~d\sigma dt \\
 &+ s  \int_{\Omega_T} (\nabla \cdot \lambda)( \nabla \cdot (\varepsilon \bar{E}))~ dxdt,~~\forall \bar{E} \in H_{E}^1(\Omega_T).\\
\end{split}
\end{equation}
Further, we obtain two equations that express that the gradients
with respect to  $\varepsilon$ and $\mu$ vanish:
\begin{equation} \label{grad1} 
\begin{split}
0 = \frac{\partial L}{\partial \varepsilon}(u)(\bar{\varepsilon})
 &=  -\int_{\Omega_T} \frac{\partial \lambda }{\partial t} 
\frac{\partial E}{\partial t}\bar{\varepsilon}~ dxdt 
+  s\int_{\Omega_T} (\nabla \cdot E) (\nabla \cdot \lambda) \bar{\varepsilon} ~dxdt \\
& +\gamma_1 \int_{\Omega} (\varepsilon - \varepsilon_0) \bar{\varepsilon}~dx,~ x \in \Omega,  
\end{split}
\end{equation}
\begin{equation} \label{grad2} 
0 = \frac{\partial L}{\partial \mu}(u)(\bar{\mu})
 =  -\int_{\Omega_T} \frac{1}{\mu^2}~\nabla \times E~ \nabla \times \lambda ~ \bar{\mu} ~dxdt 
 +\gamma_2 \int_{\Omega} (\mu - \mu_0)~ \bar{\mu}~dx,~ x \in \Omega.
\end{equation}
The equation (\ref{forward}) is the weak formulation of the state equation
(\ref{E_gauge}) and the equation (\ref{control}) is the weak
formulation of the following adjoint problem
\begin{equation}
\begin{split} \label{adjoint}
\varepsilon \frac{\partial^2 \lambda}{\partial t^2} + 
  \nabla \times (\mu^{-1} \nabla \times \lambda) -  s \varepsilon \nabla  (
 \nabla \cdot \lambda) &= -  (E - \tilde{E}) z_{\delta}, ~  x \in S_T,   \\
\lambda(\cdot, T)& =  \frac{\partial \lambda}{\partial t}(\cdot, T) = 0, \\
\partial _{n} \lambda& = \partial _{t} \lambda,~\mbox{on}~ S_{1,2},
\\
\partial _{n} \lambda& =\partial _{t} \lambda,~\mbox{on}~ S_2, \\
\partial _{n} \lambda& =0,~\mbox{on}~ S_3.
\end{split}
\end{equation}
We note that the adjoint problem is solved backward in time.

\section{Numerical method}
\subsection{Finite element discretization}
\label{sec:fem}

 We discretize $\Omega_{\rm FEM} \times (0,T)$ denoting by $K_h = \{K\}$ a partition of
 the domain $\Omega_{\rm FEM}$ into tetrahedra $K$ ($h=h(x)$ being a mesh function,
 defined as $h |_K = h_K$, representing the local diameter of the elements),
 and we let $J_k$ be a partition of $(0,T)$ into time
 intervals $J=(t_{k-1},t_k]$ of uniform length $\tau = t_k - t_{k-1}$. We
 assume also a minimal angle condition on the $K_h$ \cite{Brenner}.

To formulate the finite element method,  we
 define the finite element spaces $V_h$, $W_h^E$ and $W_h^{\lambda}$.
First we introduce the finite element trial space $W_h^E$ for every component of the electric field $E$ defined by
\begin{equation}
W_h^E := \{ w \in H_E^1: w|_{K \times J} \in  P_1(K) \times P_1(J),  \forall K \in K_h,  \forall J \in J_{\tau} \}, \nonumber
\end{equation}
where $P_1(K)$ and $P_1(J)$ denote the set of linear functions on $K$
and $J$, respectively.
We also introduce the finite element test space  $W_h^{\lambda}$ defined by
\begin{equation}
W_h^{\lambda} := \{ w \in H_{\lambda}^1: w|_{K \times J} \in  P_1(K) \times P_1(J),  \forall K \in K_h,  \forall J \in J_{\tau} \}. \nonumber
\end{equation}

Hence, the finite element spaces $W_h^E$ and $W_h^{\lambda}$
consist of continuous piecewise linear functions in space and time,
 which satisfy certain homogeneous initial and first order absorbing  boundary conditions.

To approximate functions 
$\mu(x)$  and $\varepsilon(x)$  we will use the space of piecewise constant functions $V_{h} \subset L_2(\Omega)$, 
\begin{equation*}
V_{h}:=\{u\in L_{2}(\Omega ):u|_{K}\in P_{0}(K),\forall K\in  K_h\}, 
\end{equation*}
where $P_{0}(K)$ is the piecewise constant function on $K$.

Next, we define $U_h = W_h^E \times W_h^{\lambda} \times V_h \times
V_h$. Usually $\dim U_{h}<\infty $ and $U_{h}\subset U^{1}$ as a
set and we consider $U_{h}$ as a discrete analogue of
the space $U^{1}.$ We introduce the same norm in $U_{h}$ as the one in
$U^{0},\left\Vert \cdot \right\Vert _{U_{h}}:=\left\Vert \cdot
\right\Vert _{U^{0}}.$ This means that in finite dimensional spaces
all norms are equivalent and in our computations we compute coefficients in the space $V_h$.
  The finite element method now reads: Find
$u_h \in U_h$, such that
\begin{equation*}
L'(u_h)(\bar{u})=0 ~\forall
\bar{u} \in U_h . 
\end{equation*}

\subsection{Fully discrete scheme}
\label{sec:discrete}

We expand $E$ and $\lambda$ in terms of the standard continuous piecewise
linear functions $\{\varphi_i(x)\}_{i=1}^M$ in space and
$\{\psi_k(t)\}_{k=1}^N$ in time
and substitute them into (\ref{E_gauge}) and  (\ref{adjoint}) to obtain the
following system of linear equations:
\begin{equation} \label{femod1}
\begin{split}
 M (\mathbf{E}^{k+1} - 2 \mathbf{E}^k  + \mathbf{E}^{k-1})  &= 
  - \tau^2  K \mathbf{E}^k - s \tau^2  C \mathbf{E}^k,   \\
M (\boldsymbol{\lambda}^{k+1} - 2 \boldsymbol{\lambda}^k + \boldsymbol{
\lambda}^{k-1}) &=  -\tau^2  S^k - \tau^2  K \boldsymbol{\lambda}^k - s \tau^2  D \boldsymbol{\lambda}^k, \\
\end{split}
\end{equation}
with initial  conditions :
\begin{eqnarray*}
E(\cdot, 0)&= \frac{\partial E}{\partial t}(\cdot, 0) = 0, \\
\lambda (\cdot,T) &= \frac{\partial \lambda}{\partial t} (\cdot,T) =0.
\end{eqnarray*}
  Here, $M$ is the block mass matrix in space, $K$ is the block stiffness
matrix corresponding to the rotation term, $C, D$ are the stiffness matrices
corresponding to the divergence terms, 
 $S^k$ is the load vector at
time level $t_k$, $\mathbf{E}^k$ and $ \boldsymbol{\lambda}^k$ denote the
nodal values of $E(\cdot,t_k)$ and $\lambda(\cdot,t_k)$, respectively, $\tau$
is the time step.

Let us define the mapping $F_K$ for the reference element $\hat{e}$
such that $F_K(\hat{e})=e$ and let $\hat{\varphi}$ be the piecewise
linear local basis function on the reference element $\hat{e}$ such
that $\varphi \circ F_K = \hat{\varphi}$.  Then the explicit formulas
for the entries in system (\ref{femod1}) at each element $e$ can be
given as:
\begin{equation*}
\begin{split}
  M_{i,j}^{e} & =    (\varepsilon ~\varphi_i \circ F_K, \varphi_j \circ F_K)_e, \\
  K_{i,j}^{e} & =   ( \mu^{-1} \nabla \times \varphi_i \circ F_K, \nabla \times \varphi_j \circ F_K)_e,\\
  C_{i,j}^{e} & =   ( \nabla\cdot (\varepsilon \varphi_i) \circ F_K, \nabla \cdot \varphi_j \circ F_K)_e,\\
 D_{i,j}^{e} & =   ( \varepsilon \nabla\cdot \varphi_i \circ F_K, \nabla \cdot \varphi_j \circ F_K)_e,\\
  S_{j,m}^{e}&= (E-\bar{E}, \varphi_j \circ F_K )_{e}, \\
\end{split}
 \end{equation*}
where $(\cdot,\cdot)_e$ denotes the $L_2(e)$ scalar product.

To obtain an explicit scheme, we approximate $M$ with the lumped mass matrix
$M^{L}$ (for further details, see \cite{Cohen}). Next, we multiply (\ref{femod1}) with
$(M^{L})^{-1}$ and get the following  explicit method:
\begin{equation}   \label{fem_maxwell}
\begin{split}
  \mathbf{E}^{k+1} = &2\mathbf{E}^k
  - \tau^2  (M^{L})^{-1} K\mathbf{E}^k \\
&- s \tau^2  (M^{L})^{-1} C \mathbf{E}^k   -\mathbf{E}^{k-1},  \\
  \boldsymbol{\lambda}^{k-1} = &-\tau^2 (M^{L})^{-1} S^k  
  + 2\boldsymbol{\lambda}^k
  - \tau^2  (M^{L})^{-1} K \boldsymbol{\lambda}^k \\
&- s \tau^2  (M^{L})^{-1} D \boldsymbol{\lambda}^k  -\boldsymbol{\lambda}^{k+1}.  \\
\end{split}
\end{equation} 

Finally, for reconstructing $\varepsilon(x)$  and $\mu(x)$ we can use
a gradient-based method with an appropriate initial guess values
$\varepsilon_0$  and $\mu_0$.   The discrete versions of the gradients with respect to
coefficients $\varepsilon$ and $\mu$ in (\ref{grad1}) and (\ref{grad2}), respectively, take the form:
\begin{equation*}
{g}_1(x) = - {\int_0}^T  \frac{\partial \lambda_h}{\partial t}
\frac{\partial E_h}{\partial t}~ dt  + s\int_{0}^T \nabla \cdot E_h \nabla \cdot \lambda_h dt + \gamma_1 (\varepsilon_h - \varepsilon_0)
\end{equation*}
and
\begin{equation*} 
{g}_2(x) = - {\int_0}^T \frac{1}{{\mu_h}^2}~\nabla \times E_h ~ \nabla \times \lambda_h~dt
+ \gamma_2 (\mu_h- \mu_0).
\end{equation*}
Here, $\lambda_h$  and $E_h$ are computed values of the adjoint
and forward problems  using explicit scheme
(\ref{fem_maxwell}), and $\varepsilon_h, \mu_h$   are approximated values of the computed
coefficients.

\subsection{The  algorithm}
\label{subsec:ad_alg}

In this algorithm we iteratively update approximations
$\varepsilon_{h}^{m}$ and $\mu_{h}^{m}$ of the function $\varepsilon_{
  h}$ and $\mu_{ h}$, respectively, where $m$ is the
number of iteration in our optimization procedure. We denote
\begin{equation*}
\begin{split}
{g}_1^m(x) = - {\int_0}^T  \frac{\partial \lambda_h^m}{\partial t}
\frac{\partial E_h^m}{\partial t}~ dt  +  s\int_{0}^T \nabla \cdot E_h^m \nabla \cdot \lambda_h^m  dt + \gamma_1 (\varepsilon_h^m - \varepsilon_0), \\
{g}_2^m(x) = - {\int_0}^T \frac{1}{{(\mu_h^m)}^2}~\nabla \times E_h^m ~ \nabla \times \lambda_h^m~dt
+ \gamma_2 (\mu_h^m - \mu_0), \hspace{1.8cm}
\end{split}
\end{equation*}
where functions $E_{h}\left(x,t,\varepsilon_{h}^{m}, \mu_h^m\right) ,\lambda _{h}\left(x,t,\varepsilon_{h}^{m}, \mu_h^m \right) $\ are computed by solving
the state and adjoint problems
with $\varepsilon:=\varepsilon_{h}^{m}$ and  $\mu:=\mu_{h}^{m}$.\\

\textbf{Algorithm}

\begin{itemize}
\item[Step 0.]  Choose the mesh $K_{h}$ in $\Omega$ and time partition $J$ of the time interval $\left(0,T\right) .$
Start with the initial approximations  $\varepsilon_{h}^{0}= \varepsilon_0$ and  $\mu_{h}^{0}= \mu_0$ and compute the
sequences of $\varepsilon_{h}^{m}, \mu_{h}^{m}$ via the following steps:

\item[Step 1.]  Compute solutions $E_{h}\left(x,t,\varepsilon_{h}^{m}, \mu_h^m\right) $ and $\lambda _{h}\left(x,t,\varepsilon_{h}^{m}, \mu_h^m\right) $ of state ( \ref{E_gauge}) and adjoint  (\ref{adjoint}) problems  on $K_{h}$ and $J.$

\item[Step 2.]  Update the coefficient $\varepsilon_h:=\varepsilon_{h}^{m+1}$ and $\mu_h:=\mu_{h}^{m+1}$
  on $K_{h}$ and $J$ using the conjugate gradient method
\begin{equation*}
\begin{split}
\varepsilon_h^{m+1} &=  \varepsilon_h^{m}  + \alpha_1 d_1^m(x),\\
\mu_h^{m+1} &=  \mu_h^{m}  + \alpha_2 d_2^m(x),
\end{split}
\end{equation*}
where $\alpha_i, i=1,2$, are step-sizes in the gradient update \cite{Peron} and
\begin{equation*}
\begin{split}
 d_1^m(x)&=  -g_1^m(x)  + \beta_1^m  d_1^{m-1}(x),\\
 d_2^m(x)&=  -g_2^m(x)  + \beta_2^m  d_2^{m-1}(x),
\end{split}
\end{equation*}
with
\begin{equation*}
\begin{split}
 \beta_1^m &= \frac{|| g_1^m(x)||^2}{|| g_1^{m-1}(x)||^2},\\
 \beta_2^m &=   \frac{|| g_2^m(x)||^2}{|| g_2^{m-1}(x)||^2},
\end{split}
\end{equation*}
where $d_1^0(x)= -g_1^0(x), d_2^0(x)= -g_2^0(x)$.
\item[Step 3.]  Stop computing $\varepsilon_{h}^{m}$ and obtain the
  function $\varepsilon_h$ if either $||g_1^{m}||_{L_{2}( \Omega)}\leq
  \theta$ or norms $||\varepsilon_{h}^{m}||_{L_{2}(\Omega)}$ are
  stabilized. Here,  $\theta$
  is the tolerance in $m$ updates  of the gradient method.

\item[Step 4.]  Stop computing $\mu_{h}^{m}$  and obtain the
  function $\mu_h$ if either $||g_2^{m}||_{L_{2}( \Omega)}\leq
  \theta$ or norms $||\mu_{h}^{m}||_{L_{2}(\Omega)}$ are
  stabilized. Otherwise set $m:=m+1$ and go to step 1. 

\end{itemize}

\section{Numerical Studies}
\label{sec:Numer-Simul}

\begin{figure}
 \begin{picture}(200,180)
\put(55,150){\includegraphics[angle=-90,width=14.0cm]{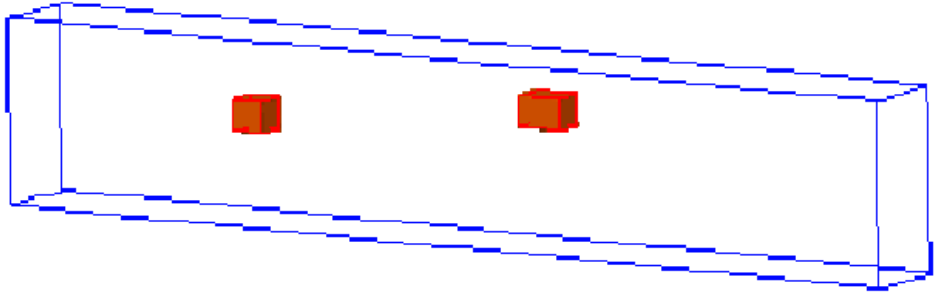}} 
 \end{picture}

 \caption{ \emph{The exact
 values of   functions $\varepsilon(x)$ and $\mu(x)$ are:
 $\varepsilon(x)=12.0, \mu(x)=2$ inside the two small scatterers, and
 $\varepsilon(x)=\mu(x)=1.0$  everywhere else in $\Omega_{\rm FEM}$.}}
 \label{fig:fig2}
 \end{figure}

In this section we present numerical simulations of the reconstruction
of two unknown functions $\varepsilon(x)$ and $\mu(x)$ inside a domain
$\Omega_{\rm FEM}$ using the algorithm of section \ref{subsec:ad_alg}.
These functions are known inside $\Omega_{\rm FDM}$ and are set to be
$\varepsilon(x) = \mu(x)=1$. 
  The goal of our numerical tests is to
reconstruct two magnetic metallic targets of figure \ref{fig:fig2}
with $\mu=2.0$.
  We note that when metallic targets are presented then
our model problem (\ref{E_gauge}) is invalid, see discussion about it
\cite{BTKF, NBKF}.  This is one of the discrepancies between our
mathematical model (\ref{E_gauge}) and the simulated backscattering
data. We refer to \cite{NBKF} for the description of other
discrepancies in a similar case.  However, one can treat metallic
targets as dielectrics with large dielectric constants and it was
shown computationally using experimental data in \cite{BTKF,
  KBKSNF,NBKF}.  Similarly with \cite{BTKF, KBKSNF, NBKF} we call
these large dielectric constants as \emph{apparent} or
\emph{effective} dielectric constants and choose values for them  in
  the interval
\begin{equation}
\varepsilon   \left(\text{metallic target}\right) \in \left(10,30\right)
.  \label{2.51}
\end{equation}
For the relative magnetic permeability we choose
values for them
  on the interval $\mu \in [1,3]$, see \cite{SSMS}.

 In our studies, we initialize only one component $E_2$ of the electrical
 field $E=(E_1,E_2,E_3)$ as the boundary condition  in (\ref{E_gauge}) on
 $S_T$ ( see (\ref{f})). Initial conditions are set to be zero.  In all 
 computations we used modification of the stabilized domain
 decomposition method of \cite{BMaxwell} which 
 was implemented using the software package WavES
\cite{waves} with two non-constant functions $\varepsilon(x)$ and $\mu(x)$.

The computational
 geometry $\Omega$ is split
 into two geometries, $\Omega_{\rm FEM}$ and $\Omega_{\rm FDM}$ such that
 $\Omega = \Omega_{\rm FEM} \cup \Omega_{\rm FDM}$, see figure \ref{fig:fig1}. Next, we introduce dimensionless spatial variables 
 $x^{\prime}= x/\left(1m\right)$ and obtain that the domain
 $\Omega_{\rm FEM}$ is transformed into  dimensionless
 computational domain
 \begin{equation*}
 \Omega_{\rm FEM} = \left\{ x= (x_1,x_2,x_3) \in (
 -3.2,3.2) \times (-0.6,0.6) \times (-0.3,0.3) \right\} .
 \end{equation*}
   The dimensionless size of our computational domain
 $\Omega$ for the forward problem is
 \begin{equation*}
 \Omega = \left\{ x= (x_1,x_2,x_3) \in (
 -3.4,3.4) \times (-0.8,0.8) \times (-0.4,0.4) \right\} .
 \end{equation*}
 The space mesh in $\Omega_{\rm FEM}$ and in $\Omega_{\rm FDM}$ consists of
 tetrahedral and cubes, respectively. 
  In the optimization algorithm  we choose the mesh size $h=0.1$
 in our geometries in the hybrid FEM/FDM method, as well as in the
 overlapping regions between FEM and FDM domains.
 In all our computational tests, we choose in (\ref{model2_1}) the penalty factor $s=1$
 in $\Omega_{\rm FEM}$.

 Note that in $\Omega_{\rm FDM}$ because of the domain decomposition
 method  and conditions (\ref{2.3}), the Maxwell's system transforms to the wave equation
 \begin{equation}\label{waveeq}
 \begin{split}
 \varepsilon \frac{\partial^2 E}{\partial t^2} - \triangle E  &=0,~ \mbox{in}~~ \Omega_T,    \\
   E(x,0) = f_0(x), ~~~E_t(x,0) &=0~ \mbox{in}~~ \Omega,     \\
  E(x,t)& = (0, f\left(t\right),0) ,~ \mbox{on}
 ~\partial \Omega_{1}\times (0,t_{1}], \\
 \partial _{n}E(x,t)& =-\partial _{t} E(x,t),~\mbox{on}
 ~\partial \Omega_{1}\times (t_{1},T), \\
 \partial _{n}E(x,t)& =-\partial _{t} E(x,t),~ \mbox{on}
 ~\partial \Omega_{2}\times (0,T), \\
 \partial _{n} E(x,t)& =0,~ \mbox{on}~\partial \Omega_{3}\times (0,T). 
 \end{split}
 \end{equation}

  We initialize only one component of the electrical field $E_2$ as a plane wave $f(t)$ 
  in
 $\Omega$ in time $T=[0,1.2]$ such that
 \begin{equation}\label{f}
 \begin{split}
 f\left(t\right) =\left\{ 
 \begin{array}{ll}
 \sin \left(\omega t \right) ,\qquad &\text{ if }t\in \left(0,\frac{2\pi }{\omega }
 \right) , \\ 
 0,&\text{ if } t>\frac{2\pi }{\omega },
 \end{array}
 \right. 
 \end{split}
 \end{equation}
 while other two components $E_1, E_3$ are initialized as zero. Thus, in $\Omega_{\rm FDM}$ we
 solve the problem (\ref{waveeq}) and in $\Omega_{\rm FEM}$ we have to solve
 \begin{equation*}
 \begin{split}
 \varepsilon \frac{\partial^2 E}{\partial t^2} + \nabla \times ( \mu^{-1} \nabla \times E)  - s\nabla  ( \nabla \cdot(\varepsilon E))  &= 0,~ \mbox{in}~~ 
 \Omega_{{\rm FEM}},    \\
   E(x,0) = 0, ~~~E_t(x,0) &= 0~ \mbox{in}~~ \Omega_{\rm FEM},    \\
 E(x,t)|_{\partial \Omega_{\rm FEM}} &= E(x,t)|_{\partial \Omega_{{\rm
       FDM}_I}}. 
 \end{split}
 \end{equation*}
Here, $\partial \Omega_{{\rm FDM}_I}$ is internal boundary of the domain $\Omega_{\rm FDM}$, and $\partial \Omega_{\rm FEM}$ is the boundary  of the domain $\Omega_{\rm FEM}$.
Similarly, in $\Omega_{\rm FDM}$ the   adjoint problem (\ref{adjoint}) transforms to  the wave equation
\begin{equation}\label{adjwaveeq}
\begin{split}
\varepsilon \frac{\partial^2 \lambda}{\partial t^2} - \triangle \lambda   &= -  (E - \tilde{E}) z_{\delta},~ \mbox{in}~~ S_T,    \\
  \lambda(x,T) = 0, ~~~\lambda_t(x,T) &=0~ \mbox{in}~~ \Omega,     \\
\partial _{n} \lambda(x,t)& =0~ \mbox{on}~S_3.
\end{split}
\end{equation}
 Thus, in $\Omega_{\rm FDM}$ we
solve the problem (\ref{adjwaveeq}) and in $\Omega_{\rm FEM}$ we have to solve
\begin{equation*}
\begin{split}
\varepsilon \frac{\partial^2 \lambda}{\partial t^2} + \nabla \times (
\mu^{-1} \nabla \times \lambda)  - s\nabla  ( \nabla \cdot(\varepsilon
\lambda))  &= 0,~ \mbox{in}~~ \Omega_{{\rm FEM}},    \\
  \lambda(x,T) = 0, ~~~\lambda_t(x,T) &= 0~ \mbox{in}~~ \Omega_{\rm FEM},    \\
\lambda(x,t)|_{\partial \Omega_{\rm FEM}} &= \lambda(x,t)|_{\partial
  \Omega_{{\rm FDM}_I}}.
\end{split}
\end{equation*}

 We  define exact functions
  ${\varepsilon(x)}=12$ and $\mu(x)=2$  inside  two small inclusions,
  see Figure \ref{fig:fig2}, and
 $\mu(x)=\varepsilon(x)=1$ at all other points of the computational domain $\Omega_{\rm FEM}$.
 We choose in our computations the time step $\tau=0.003$ which
 satisfies the CFL condition \cite{CFL67}
 and run computations in time $[0,1.2]$.

We consider the following test cases for the generation of the
backscattering data:
\begin{itemize}
\item [i)] frequency $\omega = 21$  with 3\%
  additive noise
\item [ii)] frequency $\omega = 21$  with 10\%
  additive noise
\item [iii)] frequency $\omega = 30$ with 3\%
  additive noise
\item [iv)] frequency $\omega = 30$ with 10\%
  additive noise
\end{itemize}

 To generate backscattering data at the observation points at $S_{T}$
 in each cases i)-iv), we solve the forward problem (\ref{E_gauge}),
 with function $f(t)$ given by (\ref{f}) in the time interval
 $t=[0,1.2]$ with the exact values of the parameters
 $\varepsilon(x)=12.0, \mu(x)=2$ inside scatterers of figure
 \ref{fig:fig2}, and $\varepsilon(x)=\mu(x)=1.0$ everywhere else in
 $\Omega$. We avoid the variational crime in our tests since the data
 were generated on a locally refined mesh where inclusions were
 presented.  However, our optimization algorithm works on a different 
 structured mesh with the same mesh size $h=0.1$.

\begin{figure}
  \begin{picture}(200,400)
 \put(100,-10) {\includegraphics[width=8.1cm, clip = true]{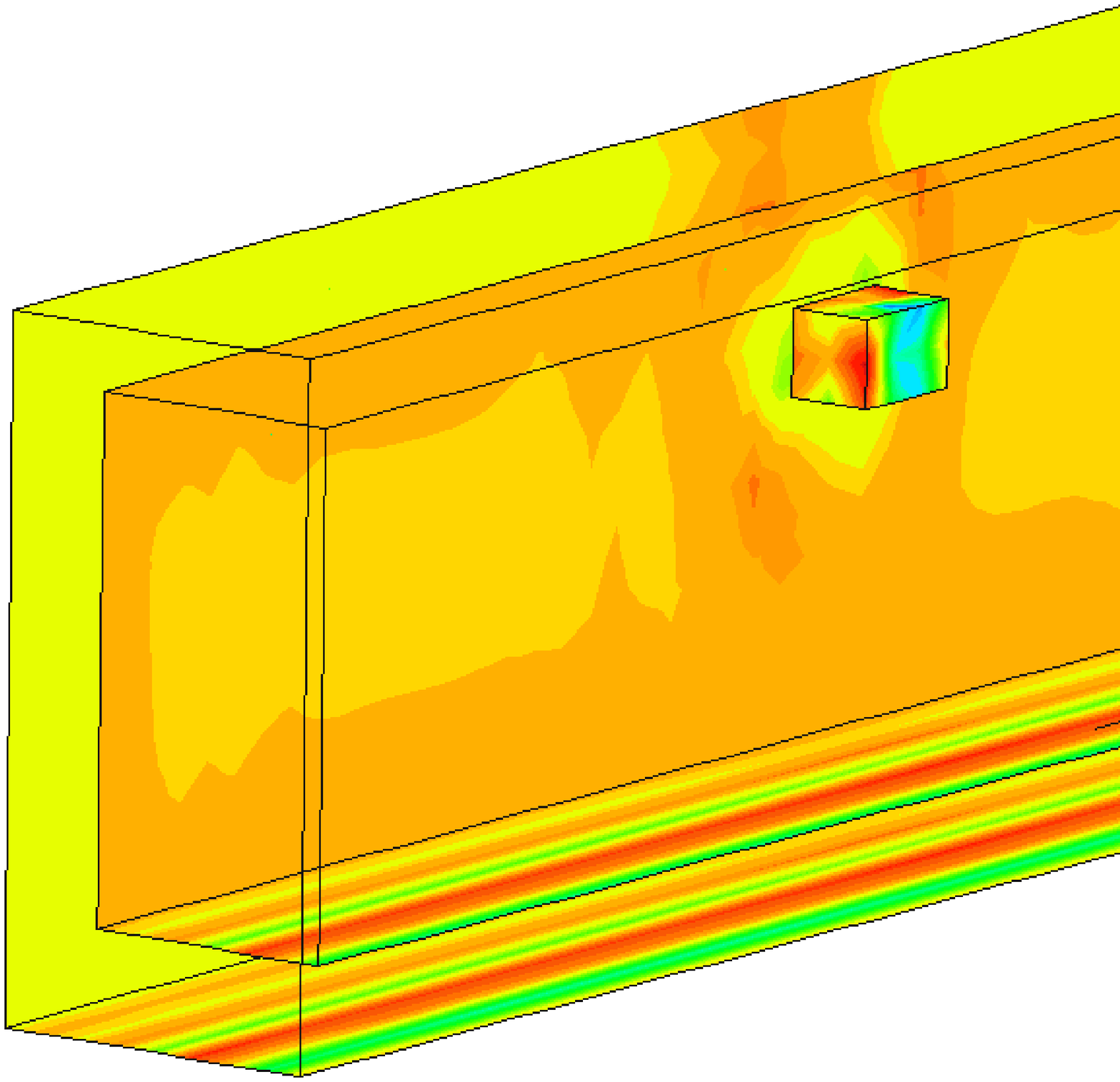}} 
 \put(330,30) {\includegraphics[width=1.5cm, clip = true]{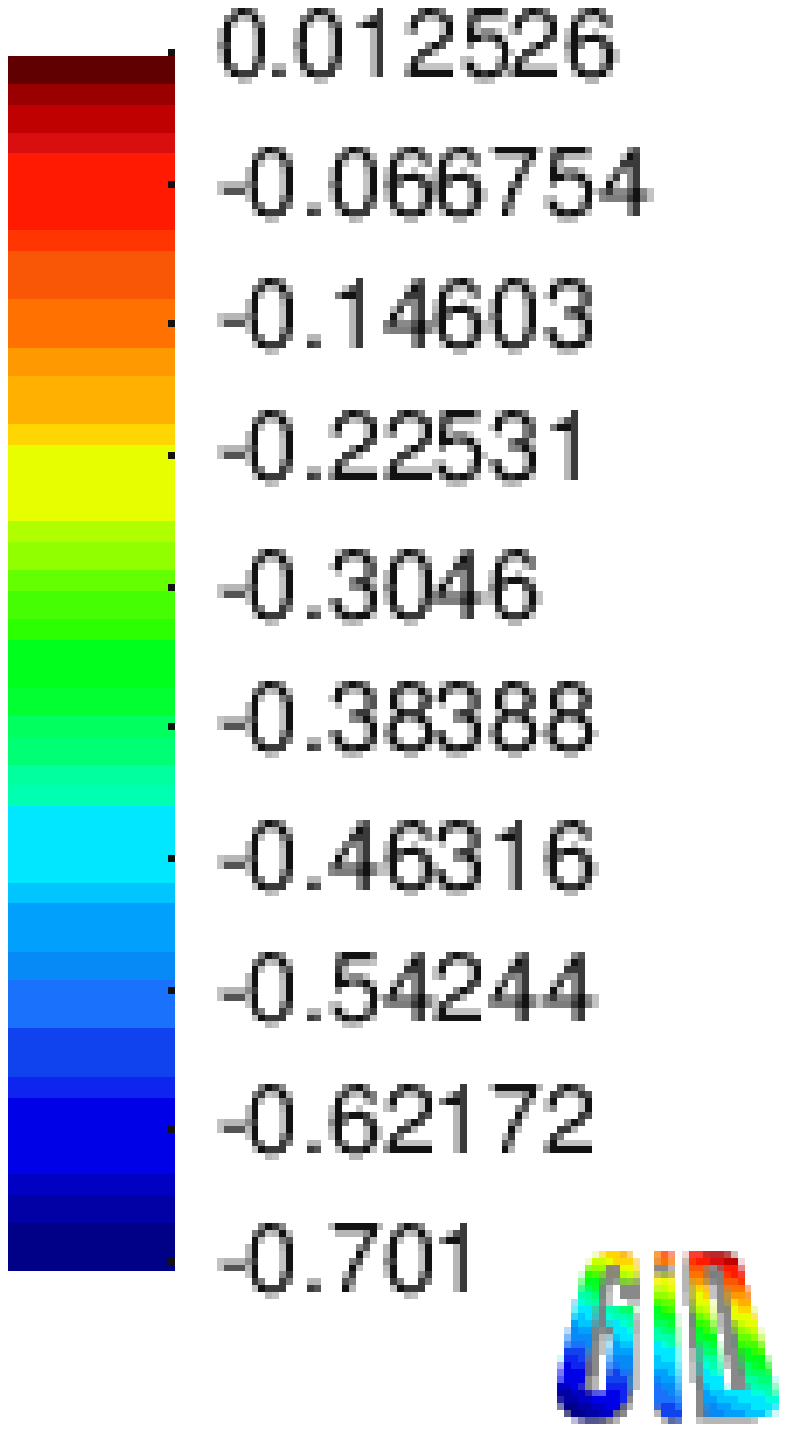}} 
 \put(100,137) {\includegraphics[width=8.1cm, clip = true]{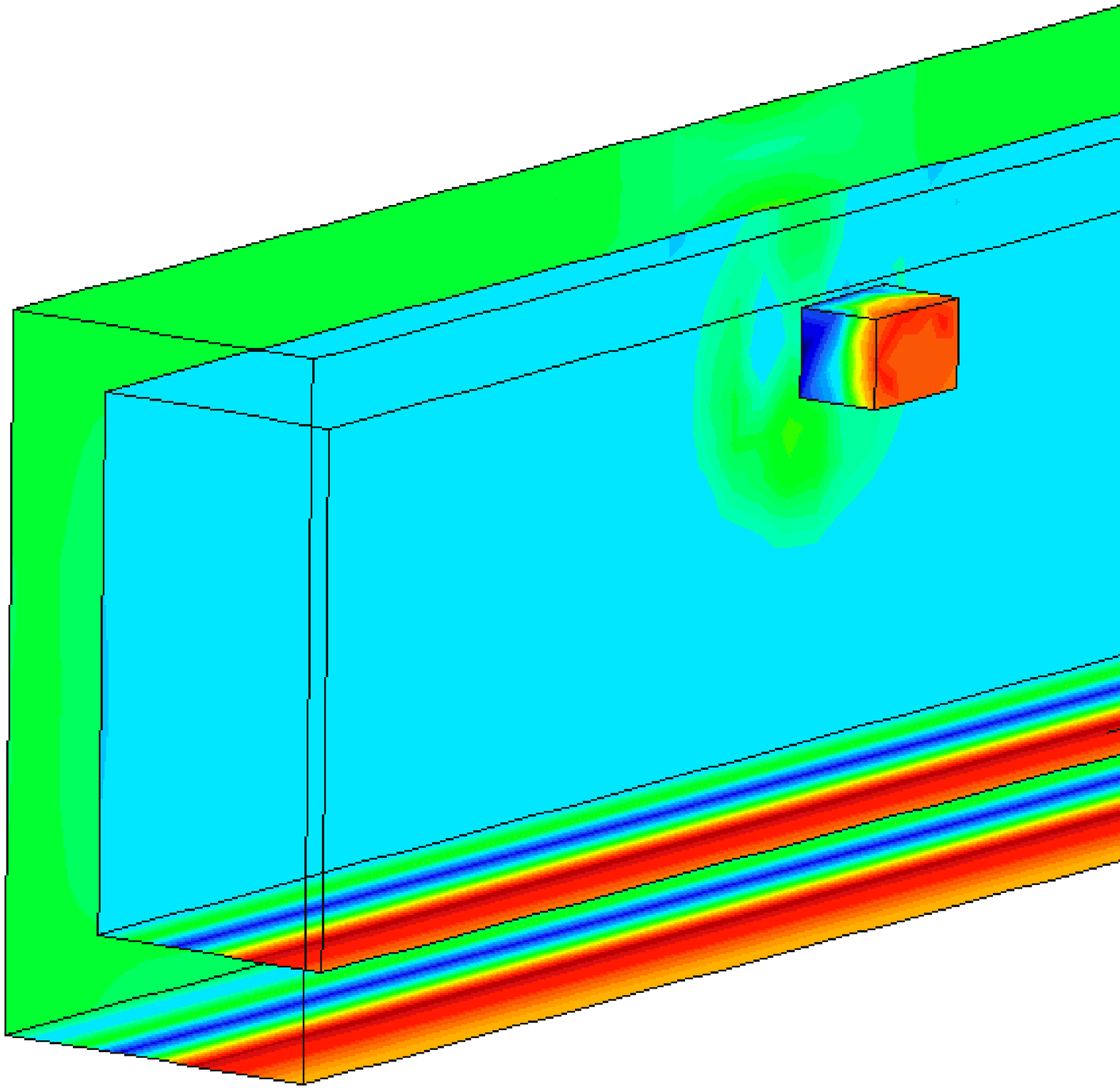}} 
 \put(330,177) {\includegraphics[width=1.5cm, clip = true]{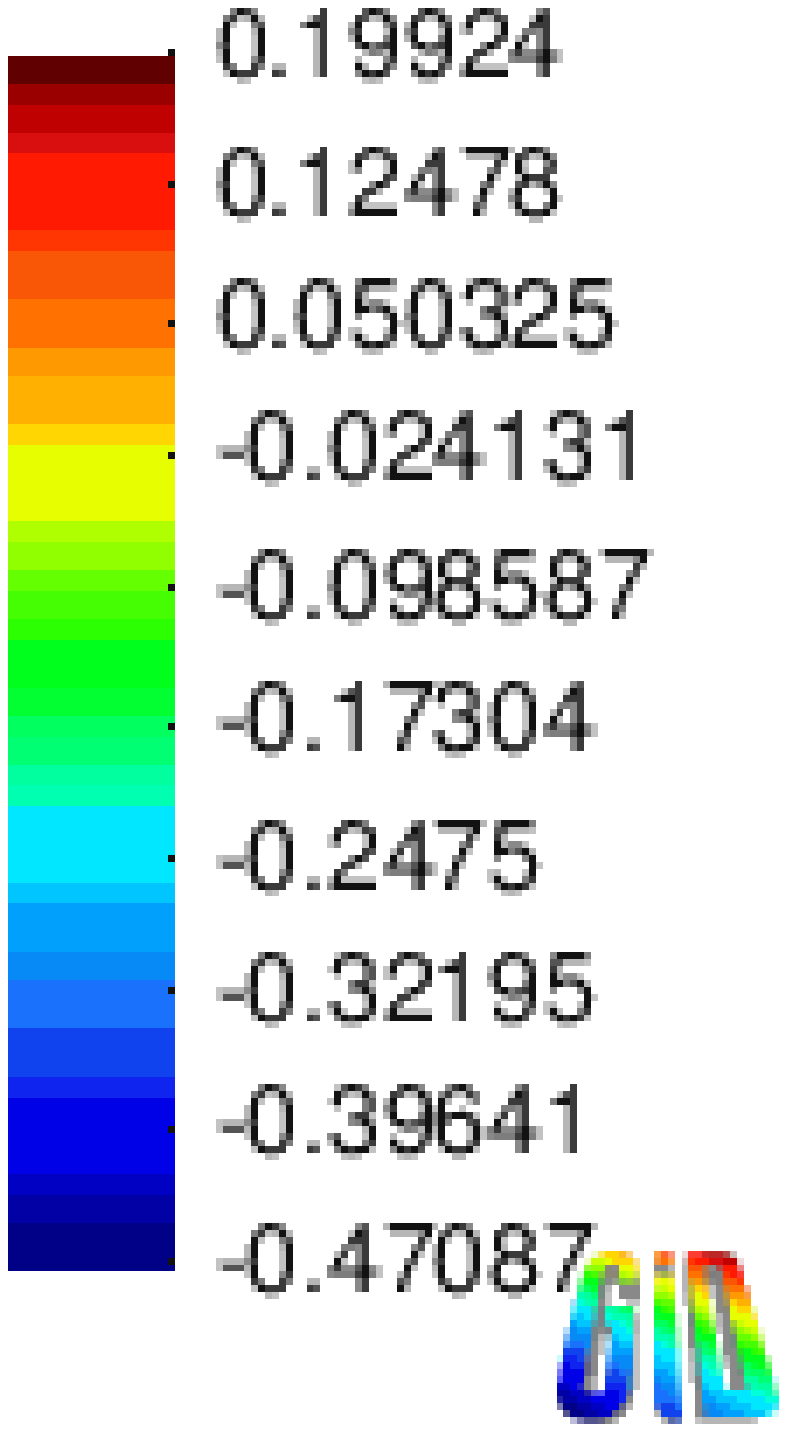}} 
 \put(100,280){\includegraphics[width=8.1cm, clip = true]{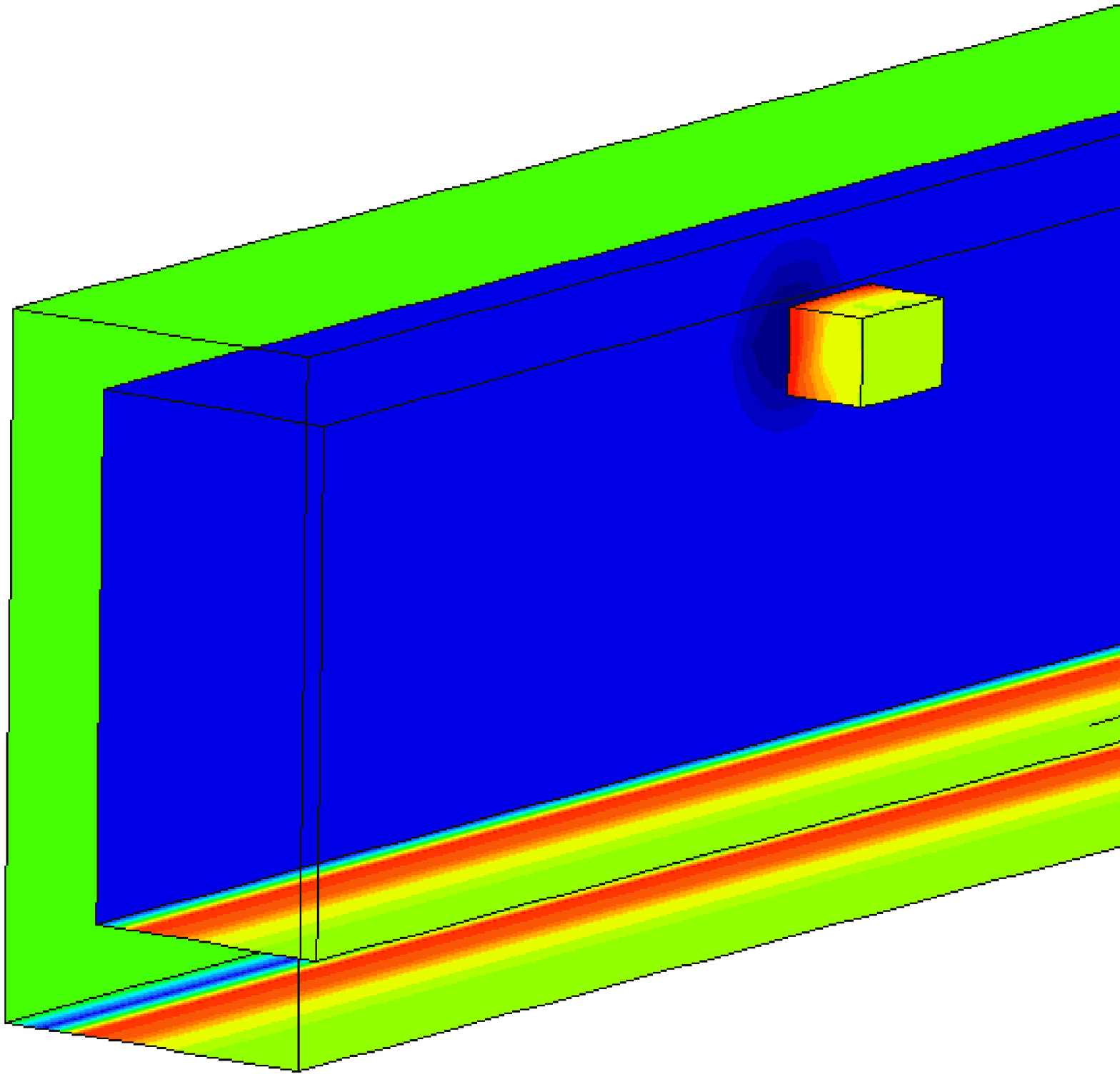}} 
 \put(330,320) {\includegraphics[width=1.5cm, clip = true]{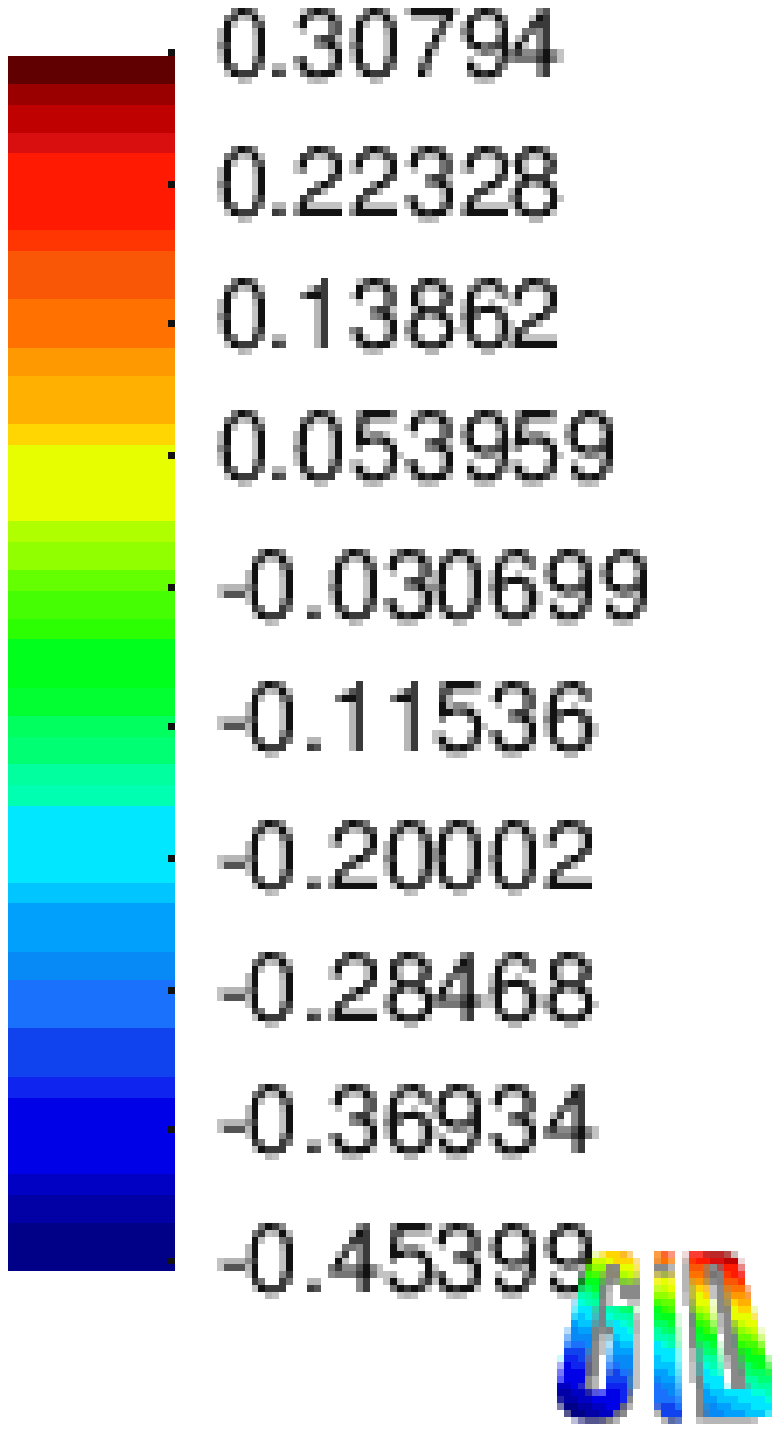}} 
\put(225,295){a) $t = 0.45$}
\put(225,150){b) $t = 0.75$}
\put(225,7){c) $t = 1.2$}
  \end{picture}
 \caption{Isosurface of the simulated exact solution of the component
   $E_2(x,t)$ for the forward problem at different times. The snapshot
   times are a) 0.45s, b) 0.75s and c) 1.2s, respectively.}
 \label{fig:Isosurfaces}
 \end{figure}

The isosurfaces of the simulated exact solution of the initialized
component $E_2(x,t)$ of the electrical field $E(x,t)$ in the forward
problem (\ref{E_gauge}) with $\omega = 30$ at different times are
presented in figure \ref{fig:Isosurfaces}. Using this figure we
observe the backscattering wave field of the component $E_2(x,t)$.

We start the optimization algorithm with guess values of the
   parameters $\varepsilon(x)=1.0$ , $ \mu(x) = 1.0$ at all points in
   $\Omega$.  Such choice of the initial guess provides a good
   reconstruction for both functions $\varepsilon(x)$ and $\mu(x)$ and
   corresponds to starting the gradient algorithm from the homogeneous
   domain, see also \cite{BKS, BMaxwell} for a similar choice of initial
   guess. Using (\ref{2.51}) the minimal and maximal values of the functions
   $\varepsilon(x)$ and $\mu(x)$ in our computations belongs to the
   following sets of admissible parameters
 \begin{equation}\label{admpar}
 \begin{split}
  M_{\varepsilon} \in \{\varepsilon\in C(\overline{\Omega })|1\leq \varepsilon(x)\leq 15\},\\
  M_{\mu} \in \{\mu\in C(\overline{\Omega })|1\leq \mu(x)\leq 3\}.
 \end{split}
 \end{equation}

The solution of the inverse problem needs to be regularized since
different coefficients can correspond to similar wave reflection data
on $\partial_1 \Omega$. We regularize the solution of the inverse
problem by starting computations with two different regularization
parameters $\gamma_1 =0.01, \gamma_2=0.9$ in (\ref{functional}).  Our
computational studies have shown that such choices for the regularization
parameters are optimal in our case.  
      We choose the regularization parameters in a computational
      efficient way such that the values of regularization parameters
      give the smallest reconstruction error given by the relative
      $L_2$ error $e_{\varepsilon} = \frac{||\varepsilon -
        \varepsilon_h ||}{||\varepsilon_h ||}$ for the reconstructed
      $\varepsilon$ and $ e_{\mu} = \frac{||\mu - \mu_h||}{||\mu_h
        ||}$ for the reconstructed $\mu$.  Here, $\varepsilon, \mu$
      are exact values of the coefficients and $\varepsilon_h, \mu_h$ are computed ones.
We refer to \cite{BKS, Engl, IJT11} for
different techniques for the choice of regularization parameters.  The
tolerance $\theta$ in our algorithm (section \ref{subsec:ad_alg}) is
set to $\theta=10^{-6}$.

\begin{figure}
 \begin{picture}(200,150)
\put(55,15){\includegraphics[width=2in]{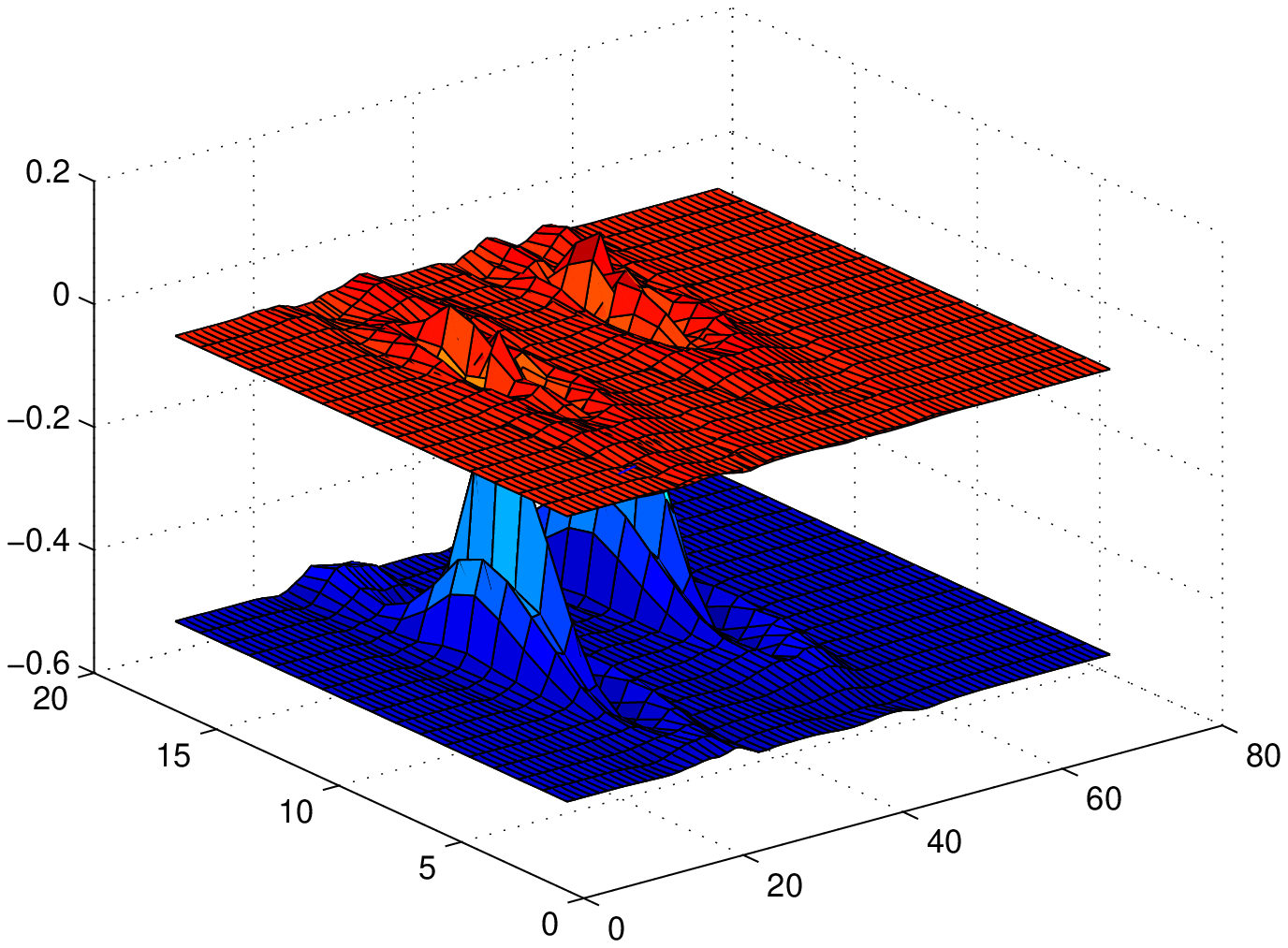}}
\put(205,15){\includegraphics[width=2in]{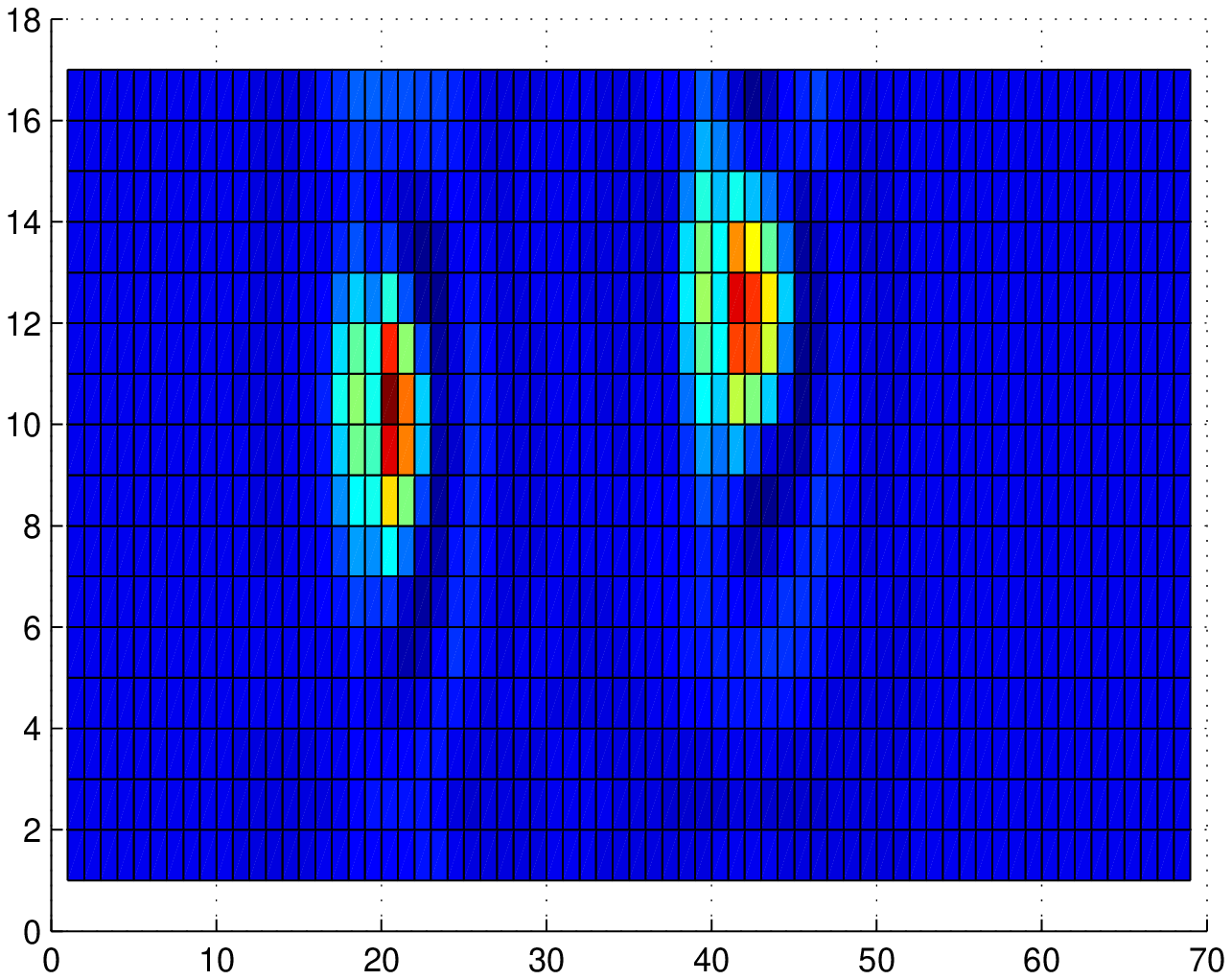}}
\put(120,-5){a)}
\put(280,-5){b)}
 \end{picture}
\caption{a) Comparison of the backscattering data of the one component
  $E_2(x,t)$ (on top) and transmitted data (below). The results are
  computed without additive noise. b) 2D projection of the
  backscattering data with 0\% additive noise. Observe that the
  backscattering data has a smaller amplitude then the transmitted
  data.}
\label{fig:omega21noise0BackScatEtTransmData}
\end{figure}

\begin{figure}
 \begin{picture}(200,300)
\put(55,15){\includegraphics[width=2in]{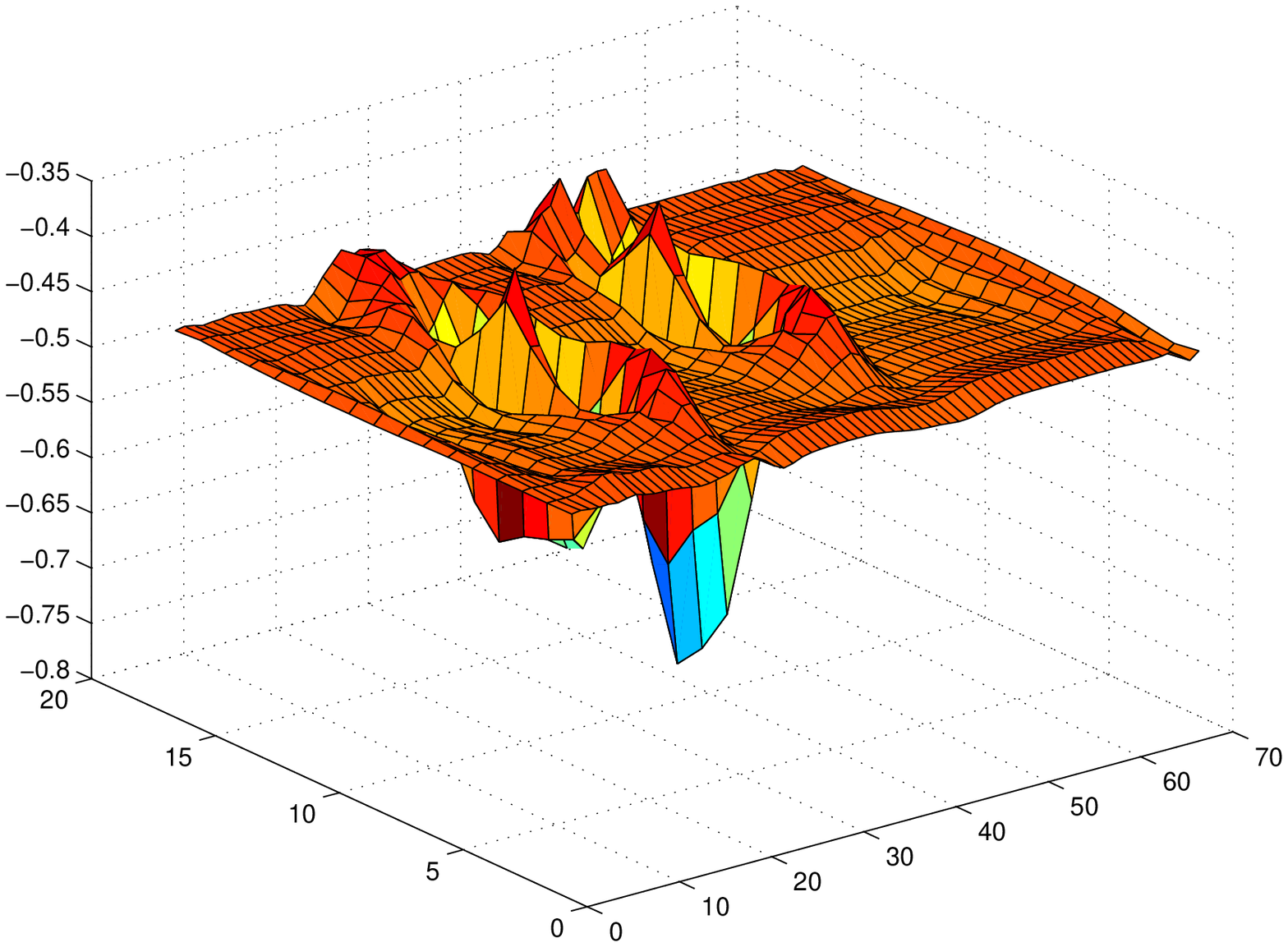}}
\put(200,15){\includegraphics[width=2.1in]{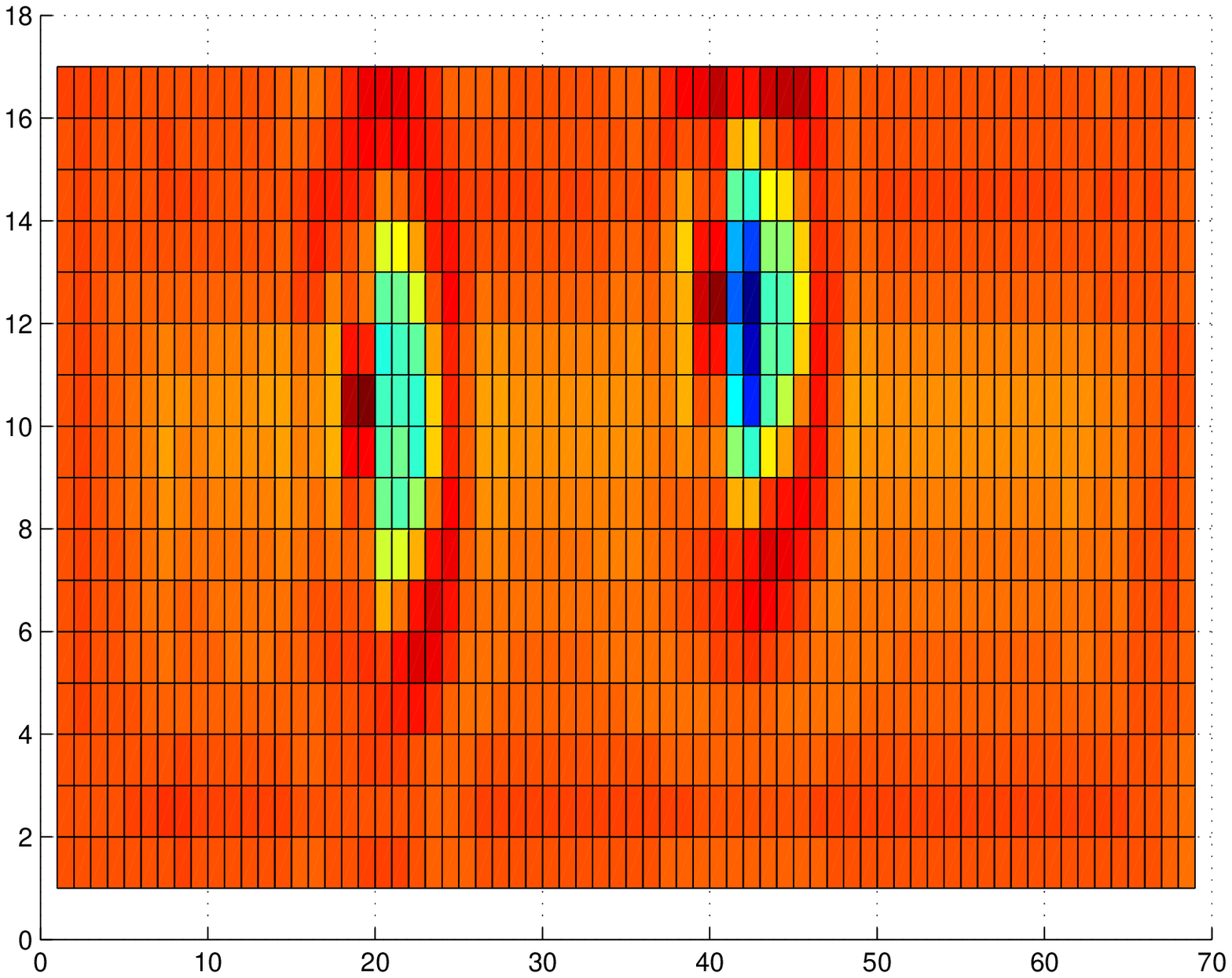}}
\put(55,175){\includegraphics[width=2in]{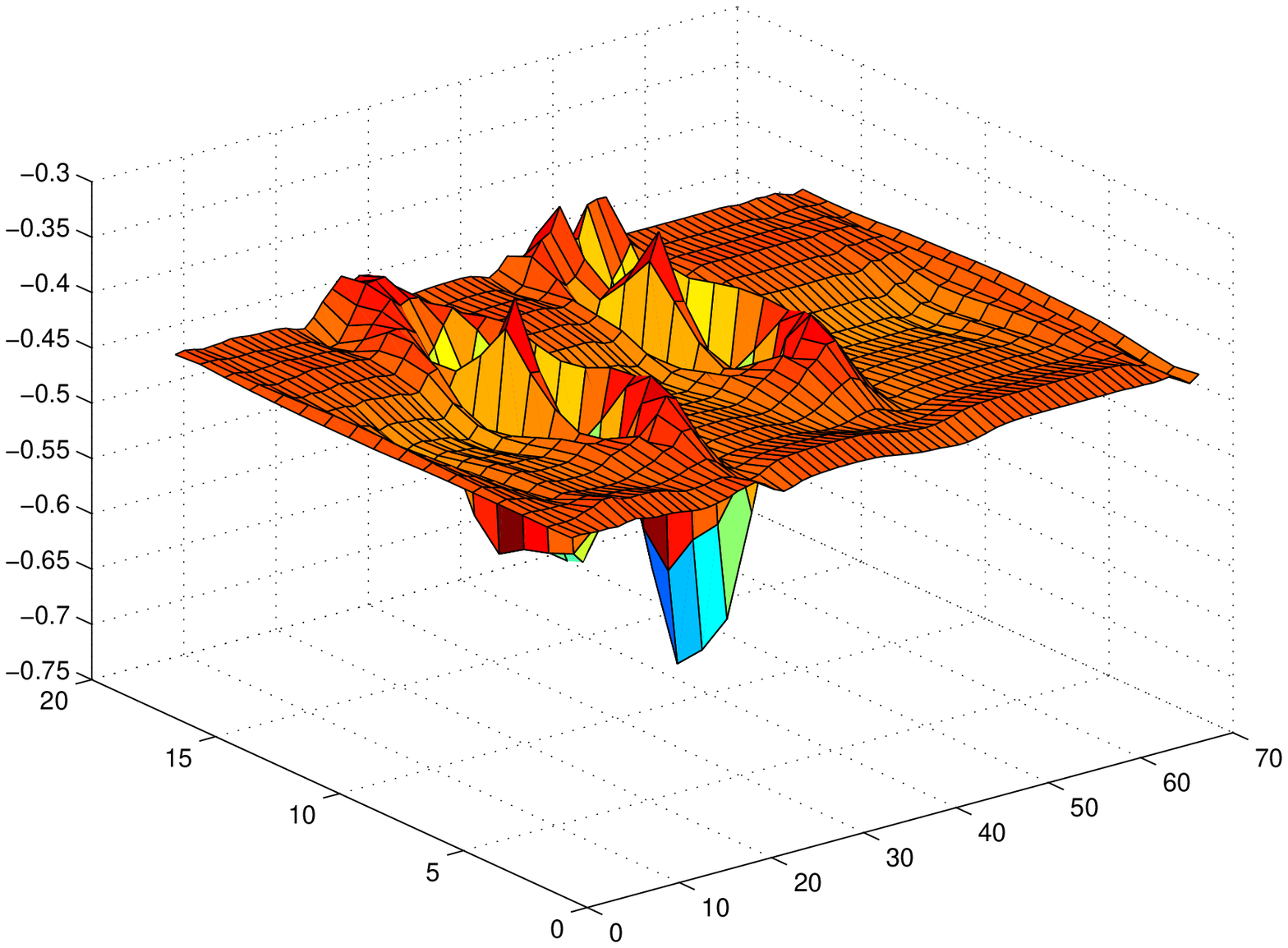}}
\put(215,175){\includegraphics[width=1.85in]{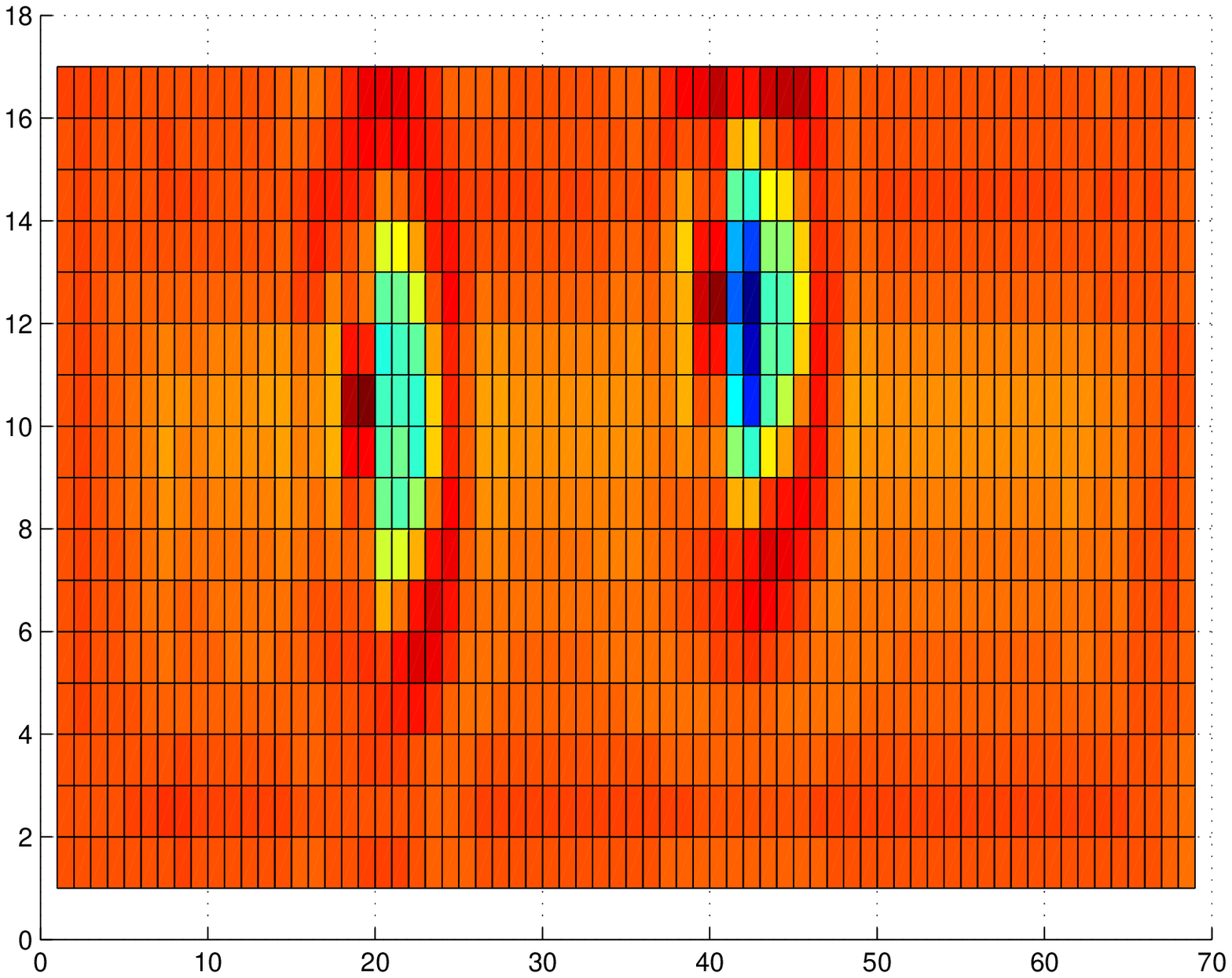}}
\put(115,-5){c)}
\put(280,-5){d)}
\put(115,160){a)}
\put(280,160){b)}
 \end{picture}
\caption{Backscattering data of the one component
  $E_2(x,t)$ using  $\omega = 21$ as the frequency. Top
  row:  3\%  additive noise  a) the   backscattering data,  b) 2D
  projection of backscattering data. Bottom row: 10\% additive noise  c) the   backscattering data,  d) 2D
  projection of backscattering data.}
\label{fig:omega21noise3et10BackScatData}
\end{figure}

\begin{figure}
 \begin{picture}(200,300)
\put(55,15){\includegraphics[width=2in]{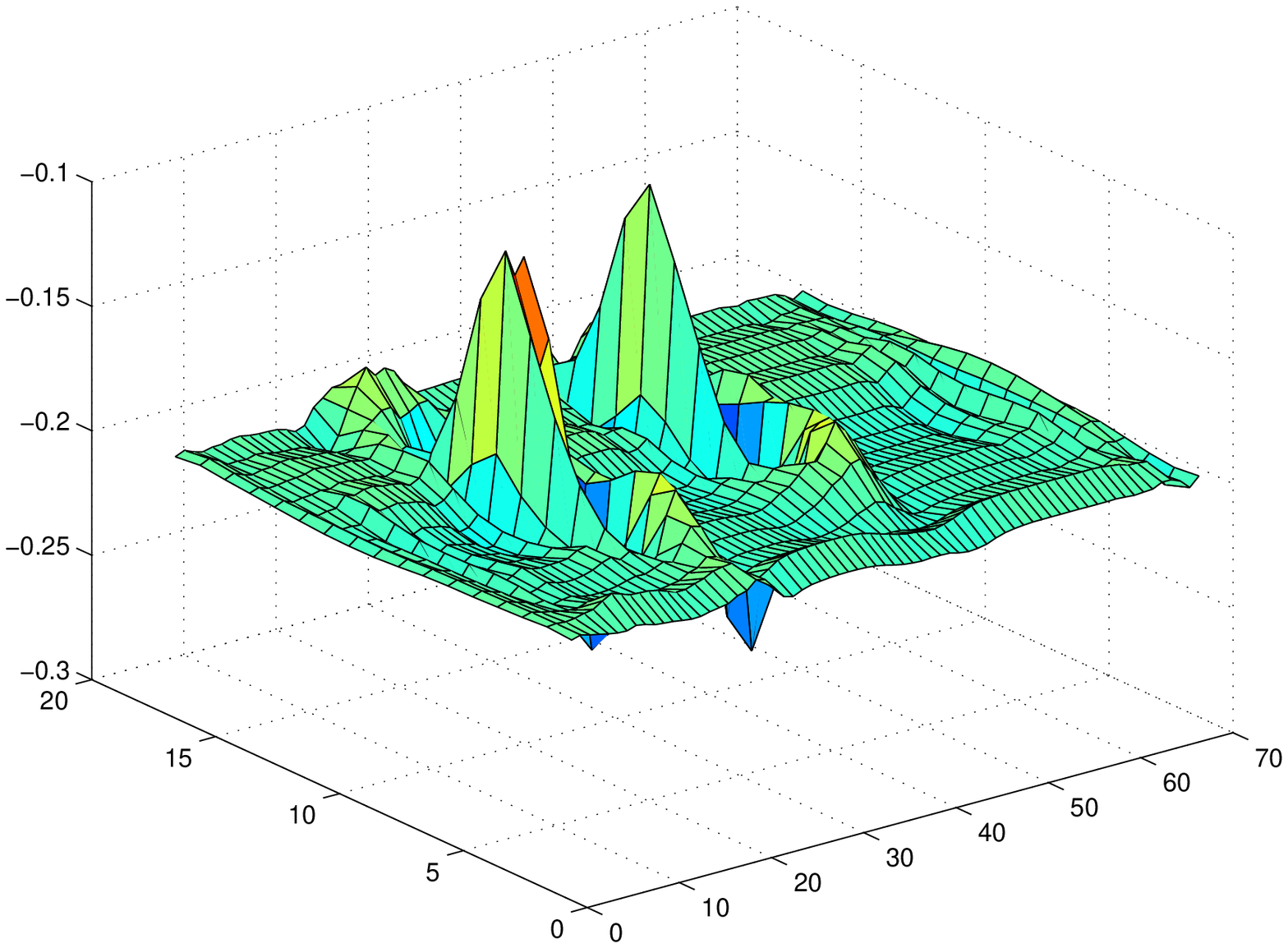}}
\put(205,15){\includegraphics[width=2in]{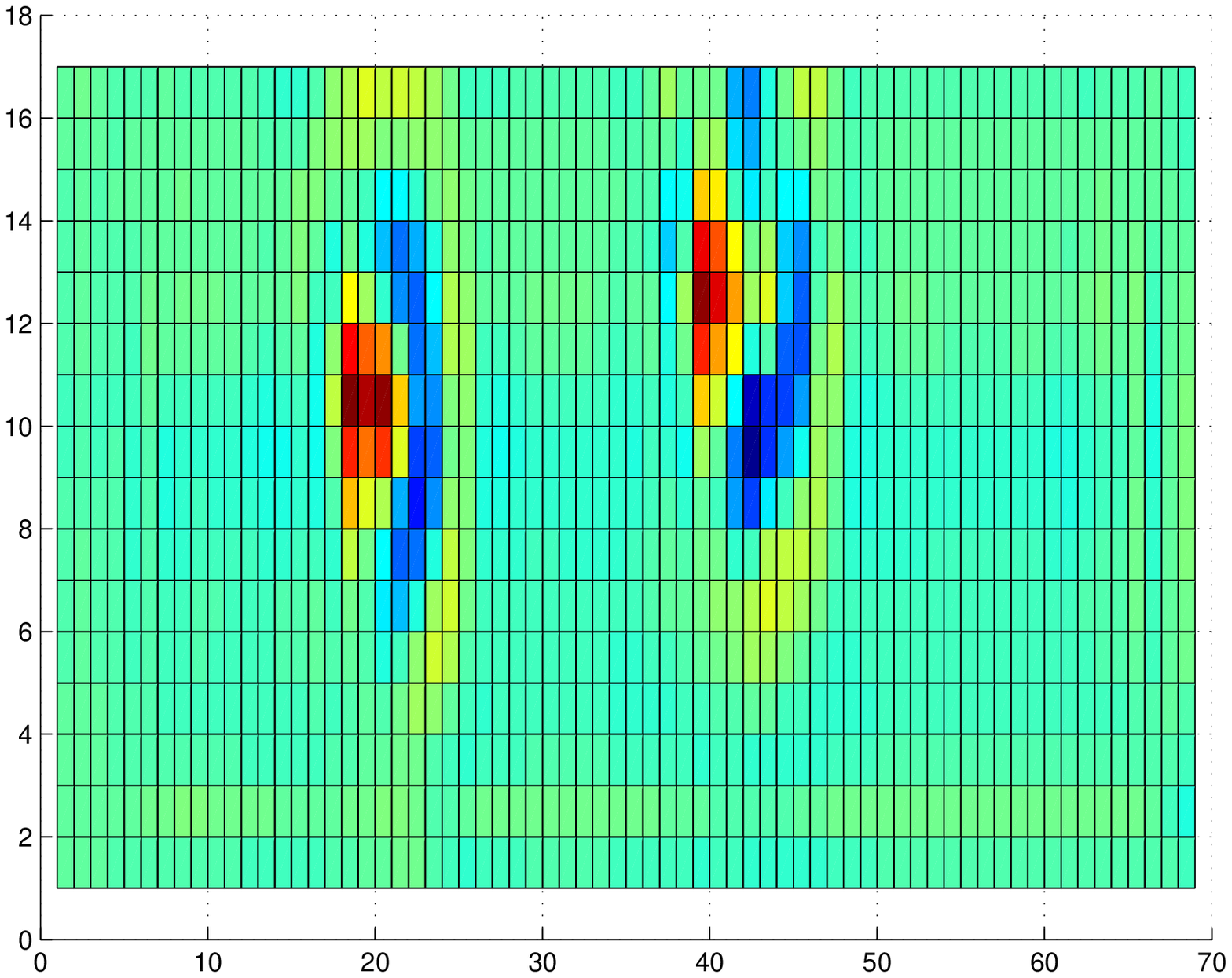}}
\put(55,175){\includegraphics[width=2in]{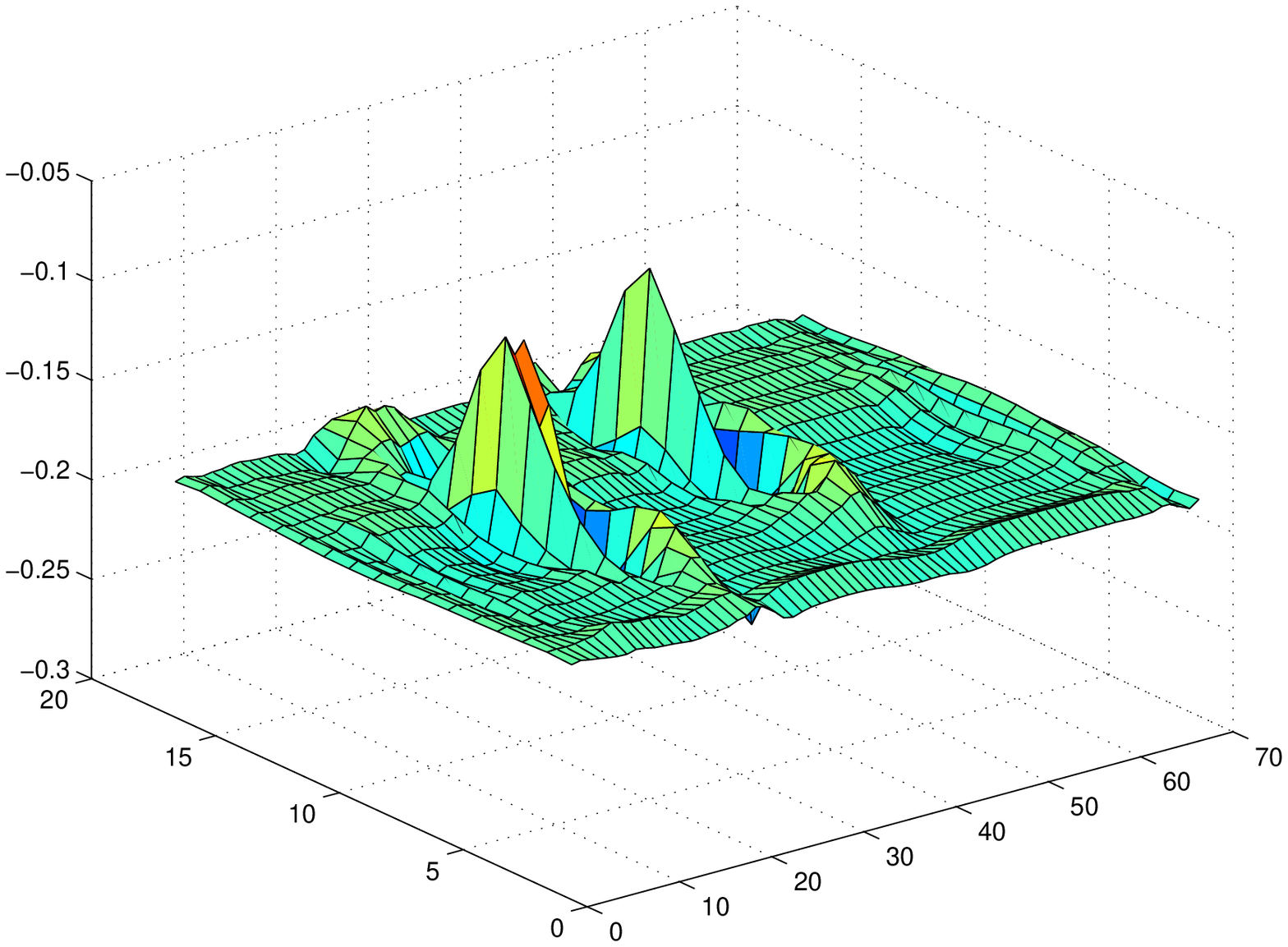}}
\put(205,175){\includegraphics[width=2in]{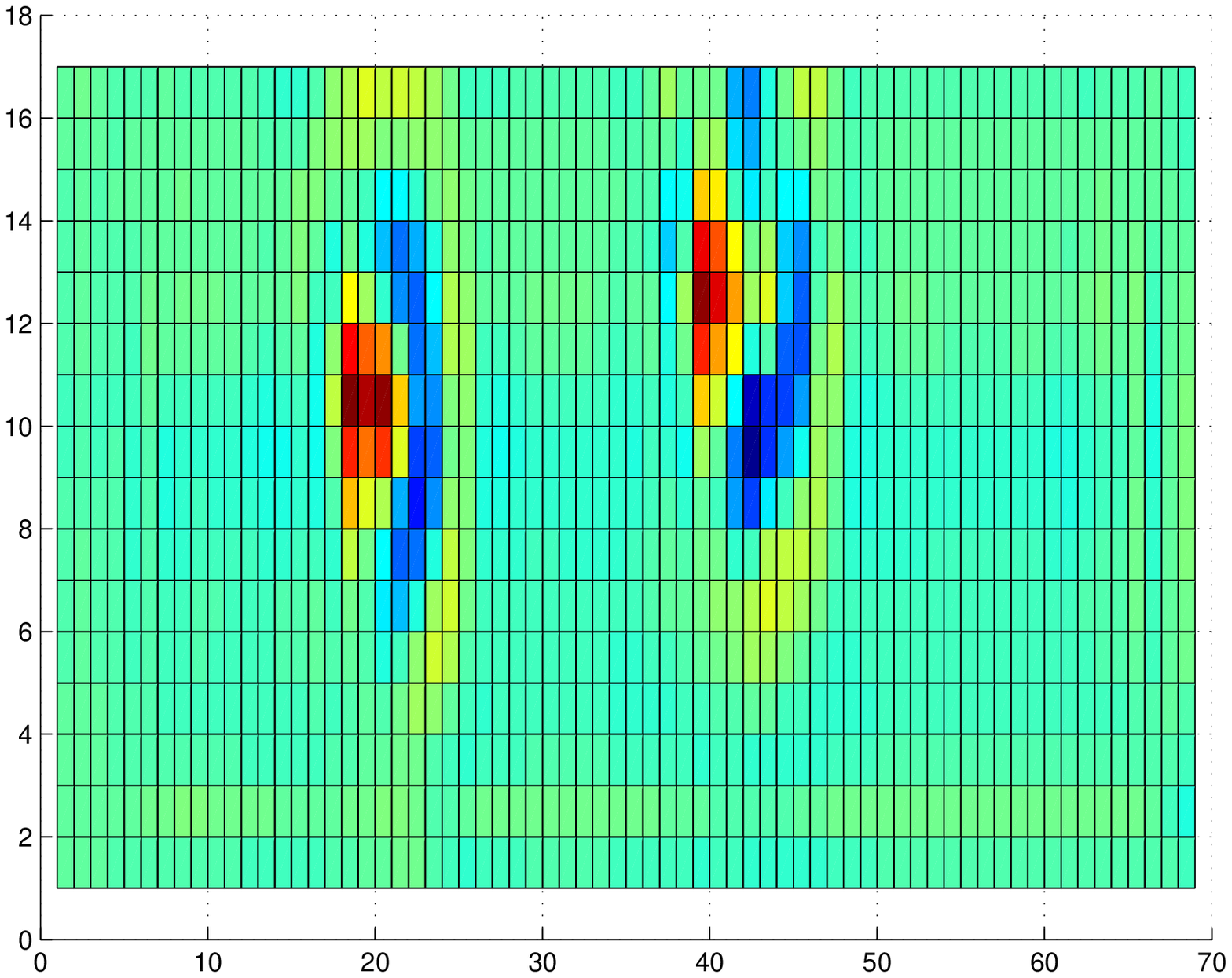}}
\put(115,-5){c)}
\put(280,-5){d)}
\put(115,160){a)}
\put(280,160){b)}

 \end{picture}
\caption{ Backscattering data of the one component
  $E_2(x,t)$ using  $\omega = 30$ as the
  frequency. Top row:  3\%
   additive noise  a) the
   backscattering data,  b) 2D projection of backscattering data. Bottom row: 10\% additive noise  c) the   backscattering data,  d) 2D
  projection of backscattering data.}
\label{fig:omega30noise3et10BackScatData}
\end{figure}

\begin{figure}
 \begin{picture}(200,180)
\put(55,15){\includegraphics[width=2in]{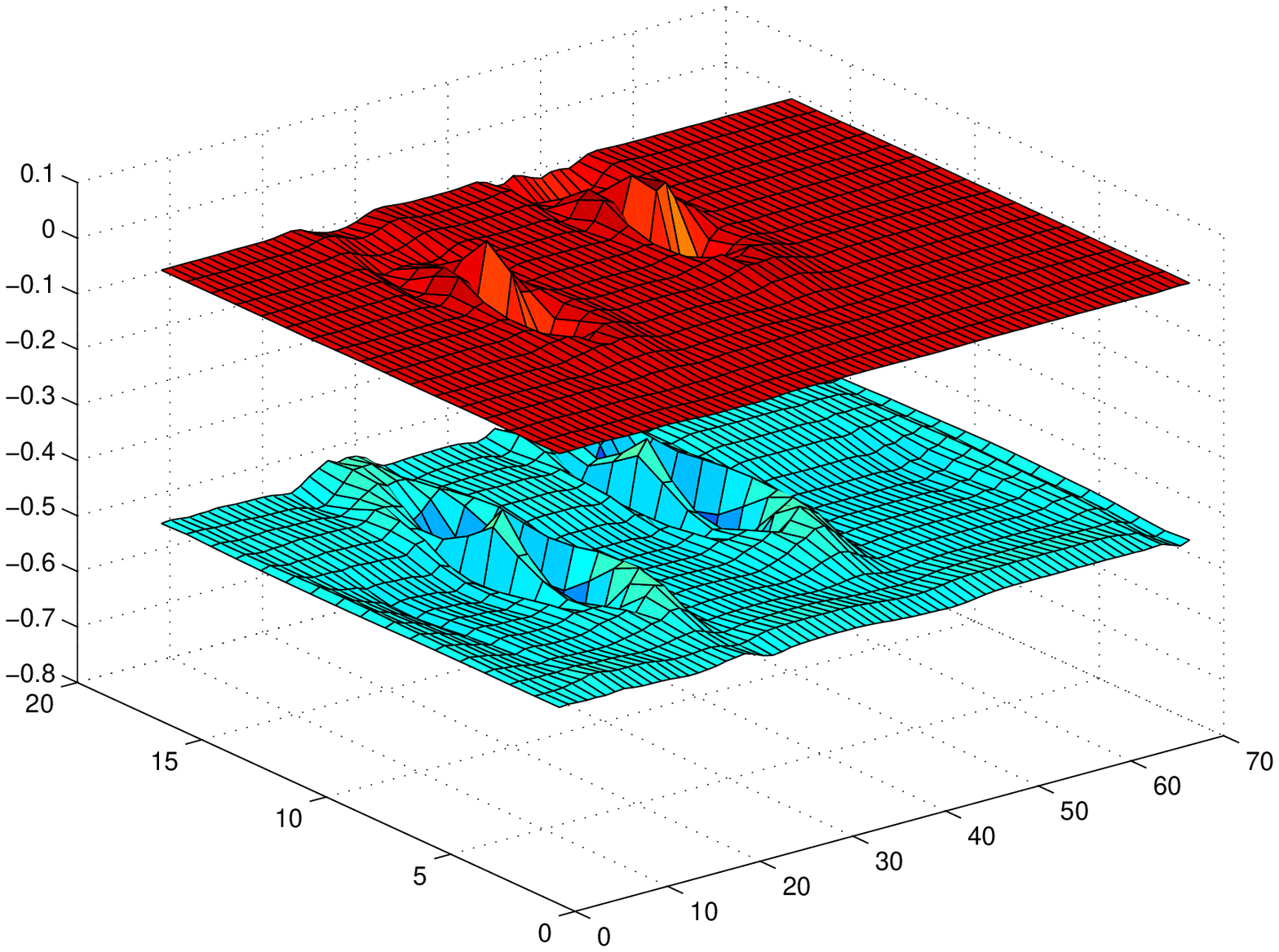}}
\put(205,15){\includegraphics[width=2in]{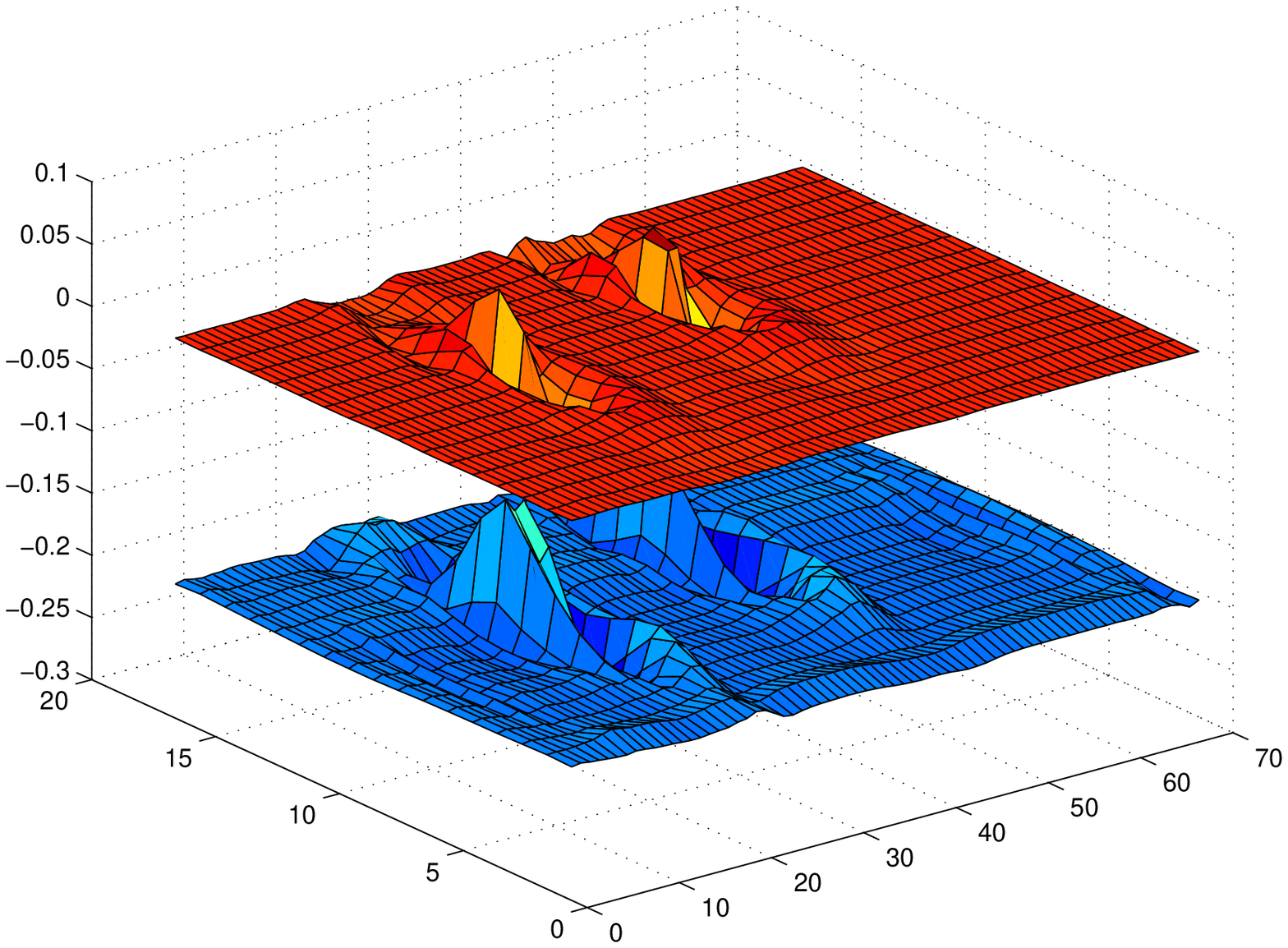}}
\put(115,-5){a)}
\put(270,-5){b)}
 \end{picture}
\caption{Comparison between computed $E_2$ (below) and $E_3$ (on top)
   components of the electric field in backscattering data with 10\%
   additive noise. a) Frequency $\omega = 21$. b) Frequency
   $\omega = 30$.}
\label{fig:CompareE2E3_omega21et30noise10}
\end{figure}

\begin{figure}
 \begin{picture}(200,180)
\put(55,15){\includegraphics[width=2in]{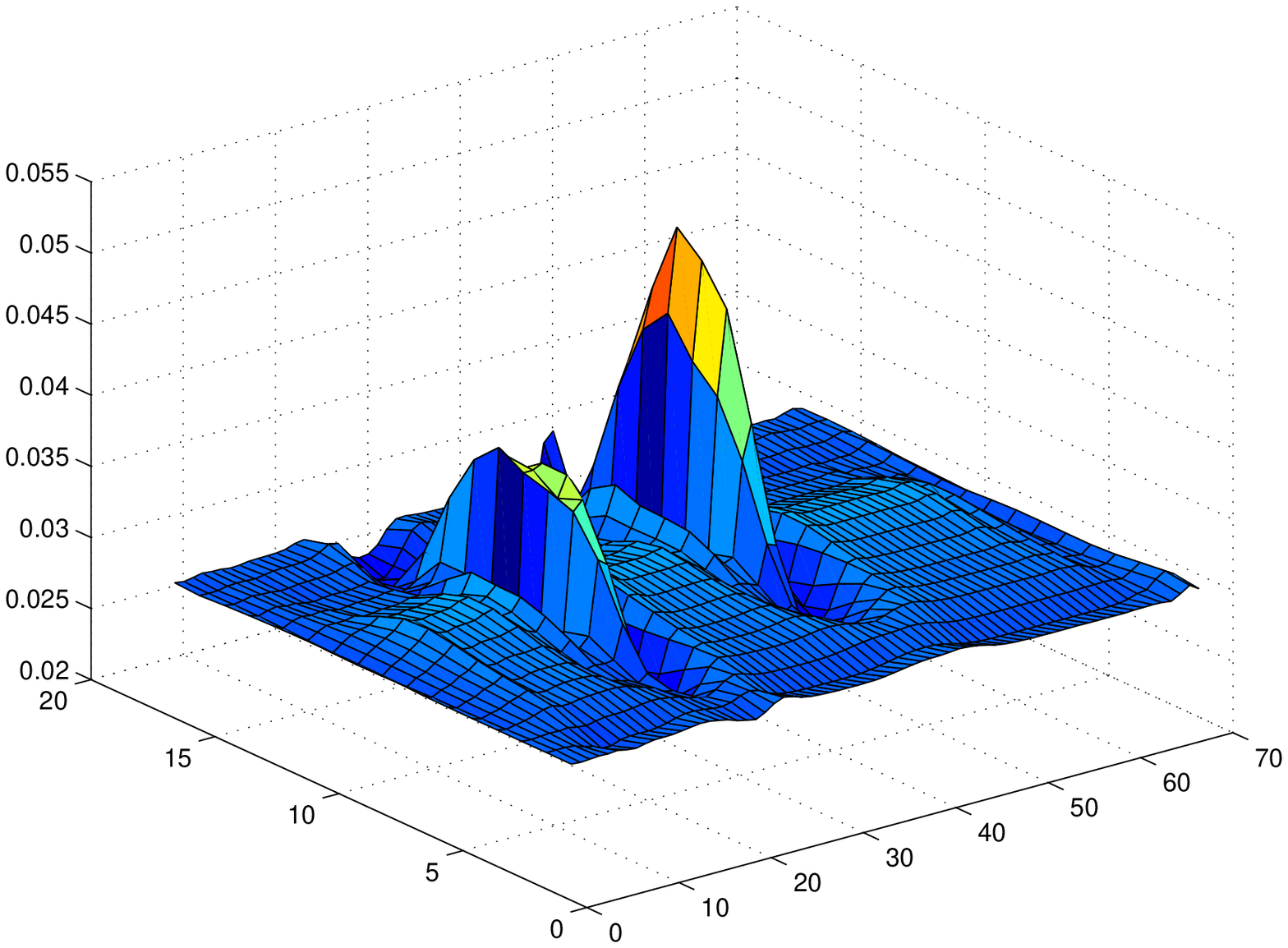}}
\put(205,15){\includegraphics[width=2in]{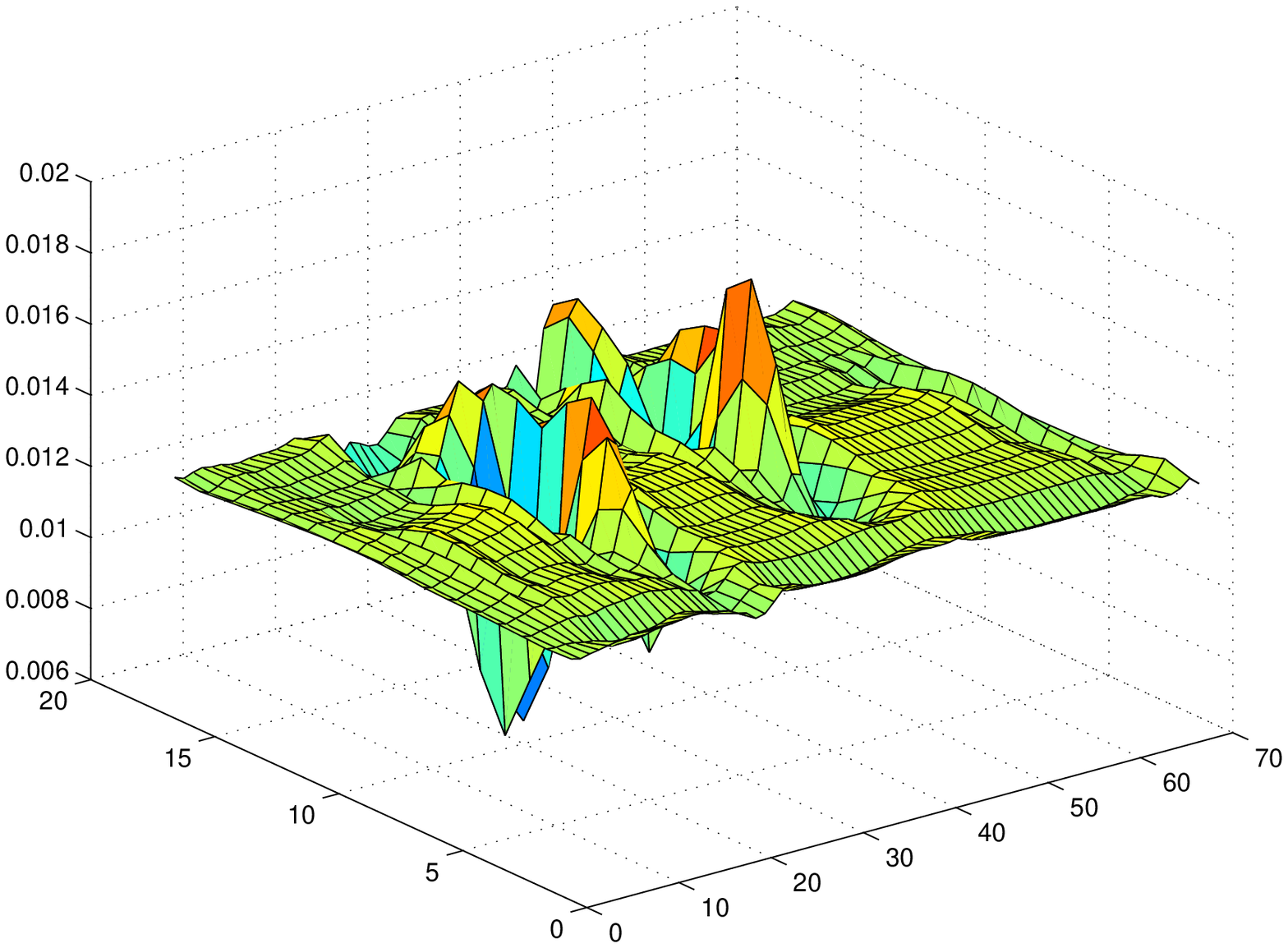}}
\put(115,-5){a)}
\put(270,-5){b)}
 \end{picture}
\caption{a) The differences of the backscattering data  of the one component
  $E_2(x,t)$ with 3\% and
   10\%
   additive noise. a) Frequency $\omega = 21$. b) Frequency $\omega = 30$.}
\label{fig:Difference_noise3noise10_omega21et30}
\end{figure}

Figure \ref{fig:omega21noise0BackScatEtTransmData}  shows a case of
backscattering data without presence of the additive noise.
 Figures
\ref{fig:omega21noise3et10BackScatData} and
\ref{fig:omega30noise3et10BackScatData} present typical behavior of
noisy backscattering data with $\omega = 21$ and $\omega = 30 $,
respectively.
Figure \ref{fig:CompareE2E3_omega21et30noise10} presents a comparison
between computed components $E_2$ and $E_3$  of the
backscattering data with 10 \% additive noise for both frequencies
$\omega = 21$ and $\omega = 30$. Figure
\ref{fig:Difference_noise3noise10_omega21et30} presents the
differences in backscattering data between 3\% and 10\% additive noise for both considered
frequencies, $\omega = 21§$ on the left and $\omega = 30$ on the right
in figure
\ref{fig:Difference_noise3noise10_omega21et30}.

\begin{figure}[tbp]
 \begin{center}
 \begin{tabular}{cc}
 {\includegraphics[angle=-90,width=5.5cm, clip = true]{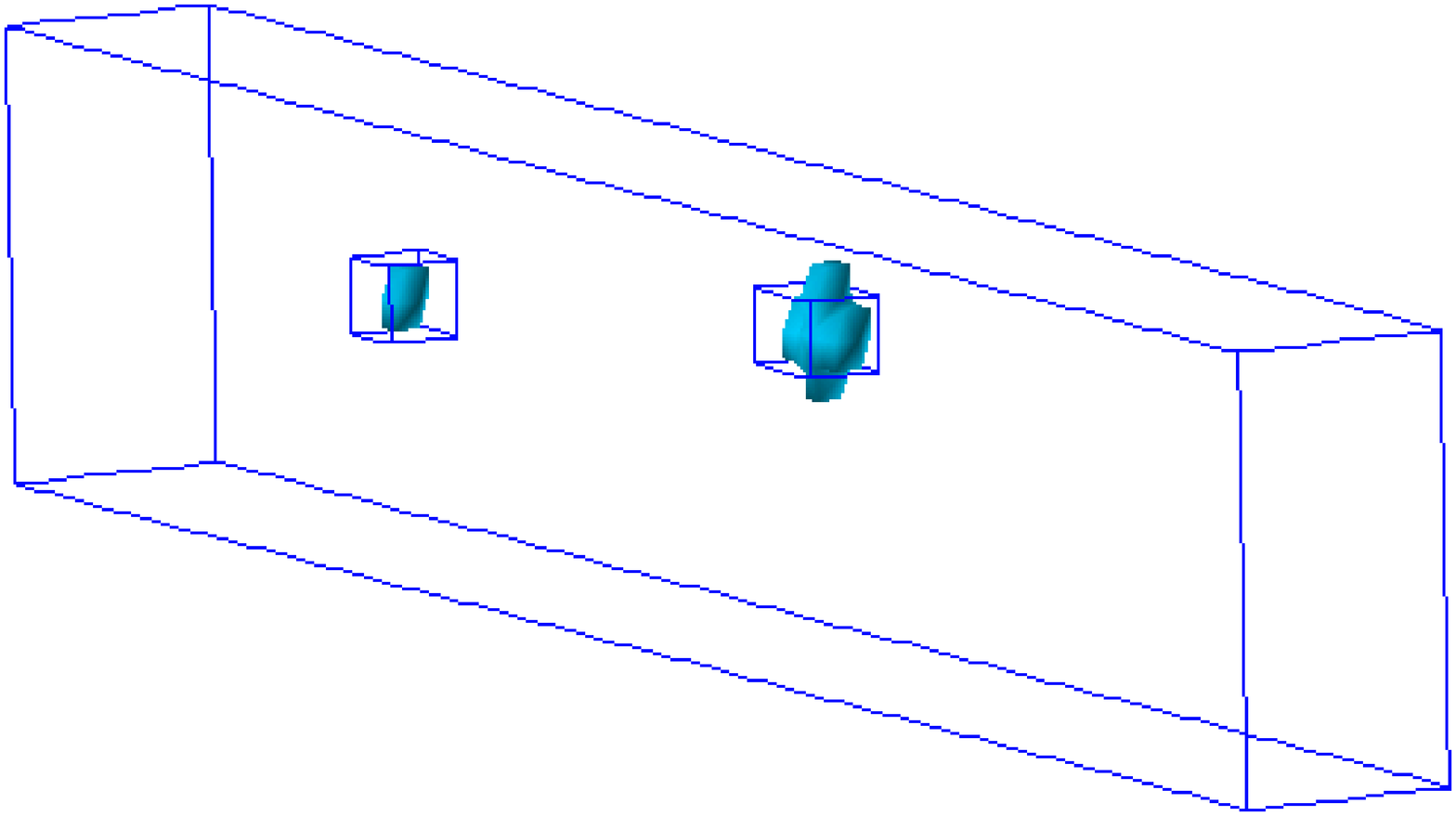}} 
 & {\includegraphics[angle=-90,width=5.5cm, clip = true]{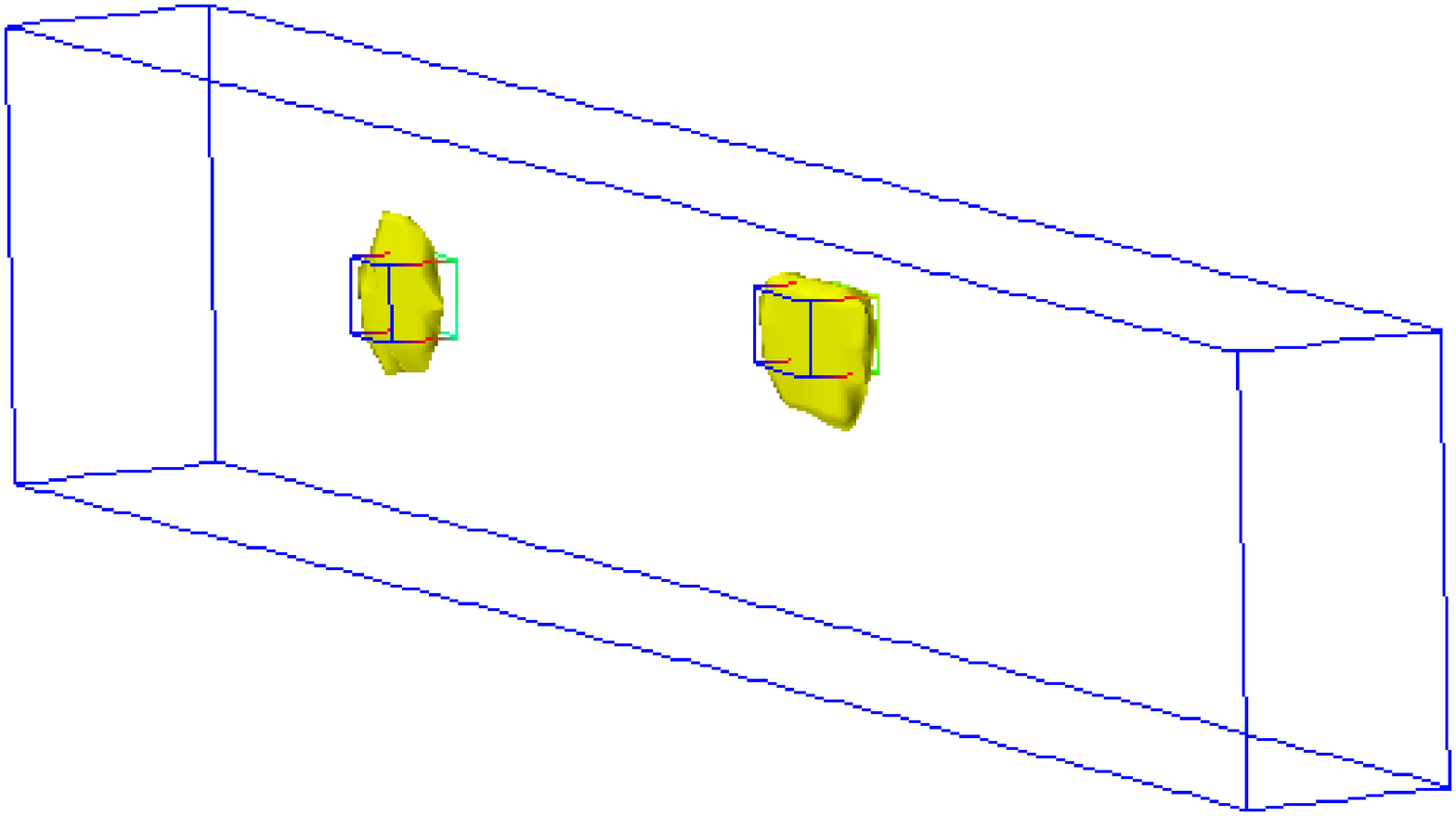}} \\
 a) $\varepsilon \approx 15$  &
 b) $\mu \approx 2$
 \end{tabular}
 \end{center}
 \caption{Computed images of  reconstructed  functions
   $\varepsilon(x)$ and $\mu(x)$ from backscattering data with 3\%
   additive noise, frequency $\omega = 21$.}
 \label{fig:Recos_omega21noise3}
 \end{figure}

The reconstructions of $\varepsilon(x)$ and $\mu(x)$ with $\omega =
21$ using 3\% and 10\% noise, are presented in figures
\ref{fig:Recos_omega21noise3} and
\ref{fig:Recos_omega21noise10}. Similarly, reconstructions of
$\varepsilon(x)$ and $\mu(x)$ with $\omega = 30$ using  3\% and 10\%
of additive noise, are presented in figures
\ref{fig:Recos_omega30noise3} and \ref{fig:Recos_omega30noise10}, respectively. 

\begin{figure}
 \begin{picture}(200,240)
\put(55,108) {\includegraphics[angle=-90,width=4.1cm, clip = true]{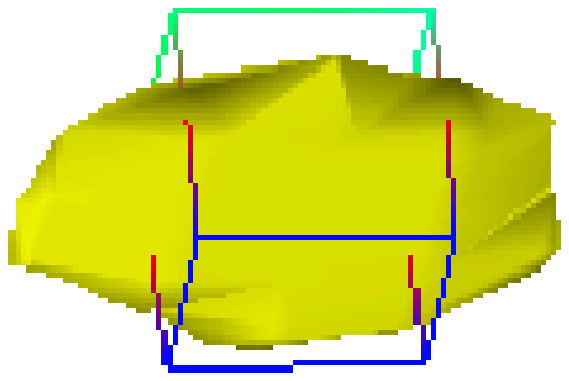}} 
\put(210,108) {\includegraphics[angle=-90,width=4.1cm, clip = true]{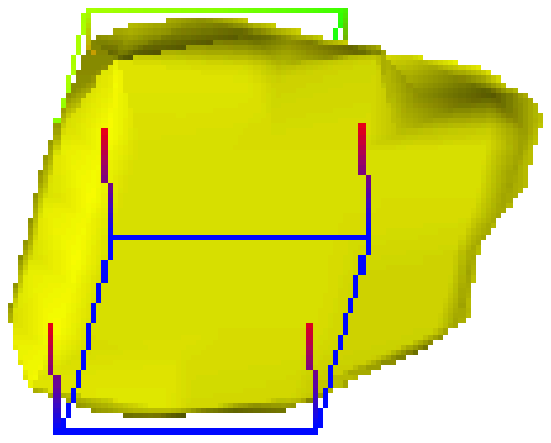}} 
\put(65,235){\includegraphics[angle=-90,width=4.5cm, clip = true]{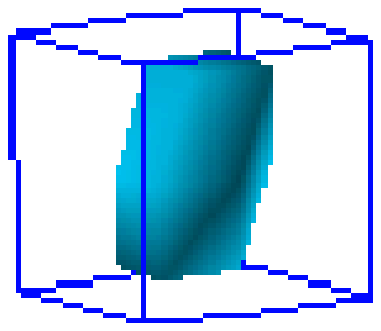}} 
\put(210,235){\includegraphics[angle=-90,width=4.5cm, clip =  true]{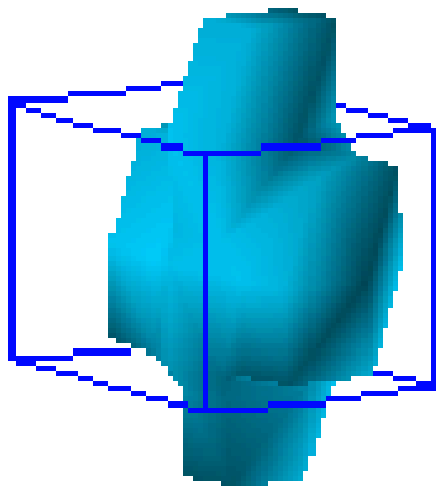}} 
 \end{picture}
\caption{Zoomed reconstructions of $\varepsilon$ (top row) and $\mu$
  (bottom row). Frequency $\omega = 21$ additive noise 3\%.  }
\label{fig:ZoomedEpsMu_omega21noise3}
\end{figure}

\begin{figure}[tbp]
 \begin{center}
 \begin{tabular}{cc}
 {\includegraphics[angle=-90,width=5.5cm, clip = true]{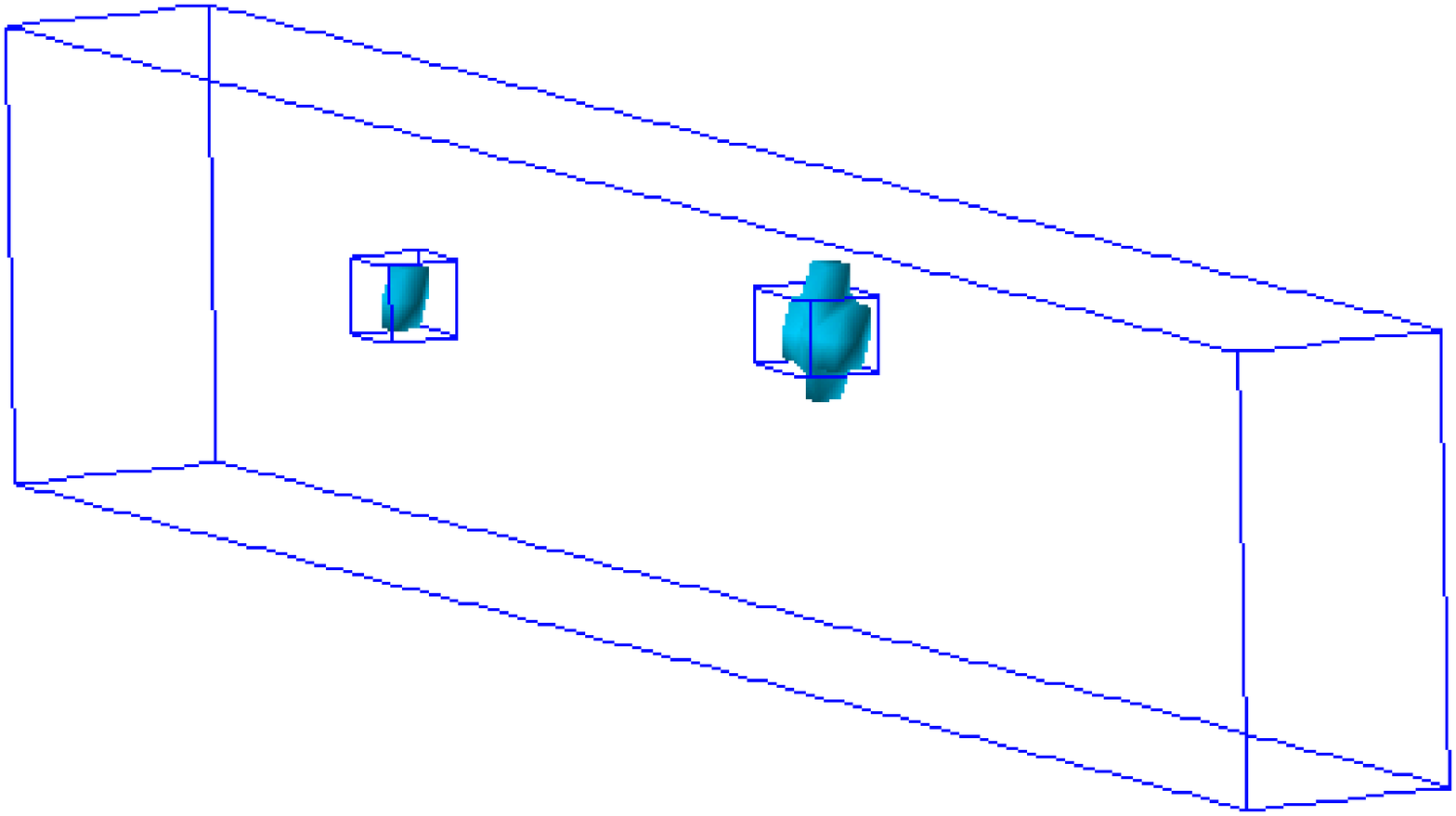}} 
 & {\includegraphics[angle=-90,width=5.5cm, clip = true]{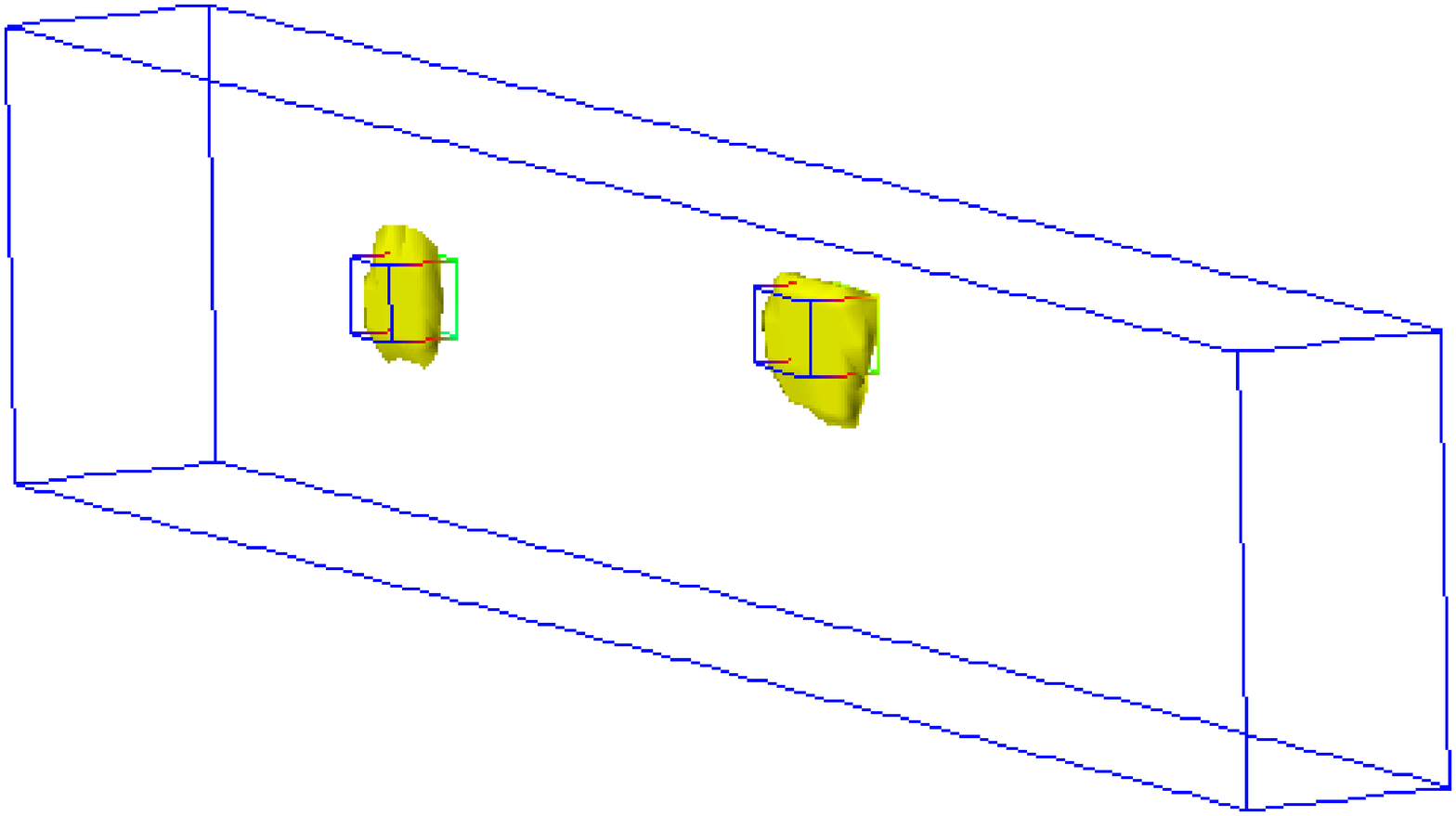}} \\
 a) $\varepsilon \approx 15$  &
 b) $\mu \approx 2$
 \end{tabular}
 \end{center}
 \caption{Computed images of  reconstructed  functions
   $\varepsilon(x)$ and $\mu(x)$ from backscattering data with 10\%
   additive noise, frequency $\omega = 21$.}
 \label{fig:Recos_omega21noise10}
 \end{figure}

\begin{figure}
 \begin{picture}(200,240)
\put(60,108) {\includegraphics[angle=-90,width=4.2cm, clip = true]{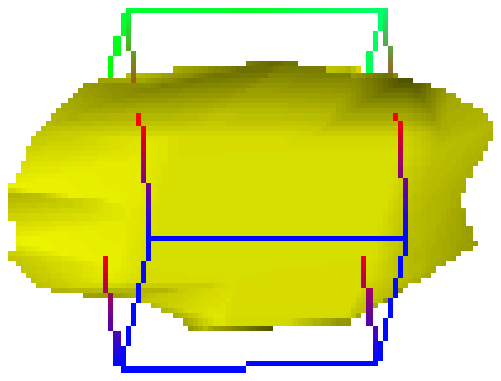}} 
\put(210,108) {\includegraphics[angle=-90,width=4.2cm, clip = true]{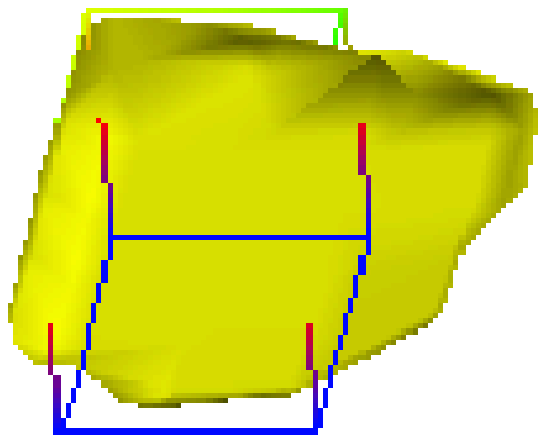}} 
\put(60,235){\includegraphics[angle=-90,width=4.5cm, clip = true]{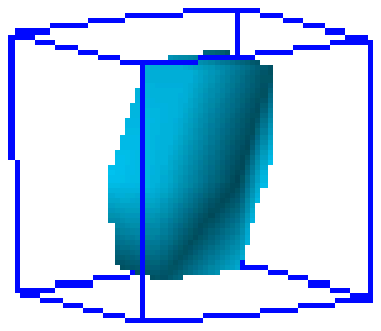}} 
\put(210,235){\includegraphics[angle=-90,width=4.5cm, clip =  true]{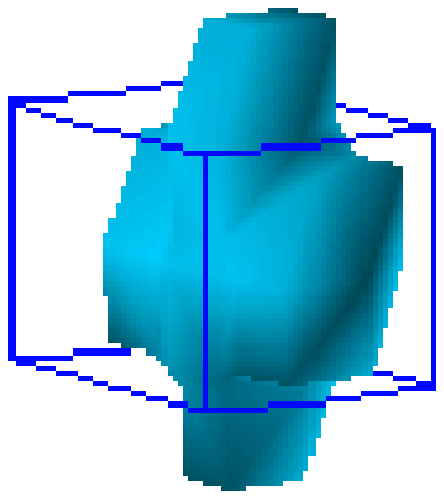}} 
 \end{picture}
\caption{Zoomed reconstructions of $\varepsilon$ (top row) and $\mu$
  (bottom row). Frequency $\omega = 21$ additive noise 10\%.  }
\label{fig:ZoomedEpsMu_omega21noise10}
\end{figure}

\begin{figure}[tbp]
 \begin{center}
 \begin{tabular}{cc}
 {\includegraphics[angle=-90,width=5.5cm, clip = true]{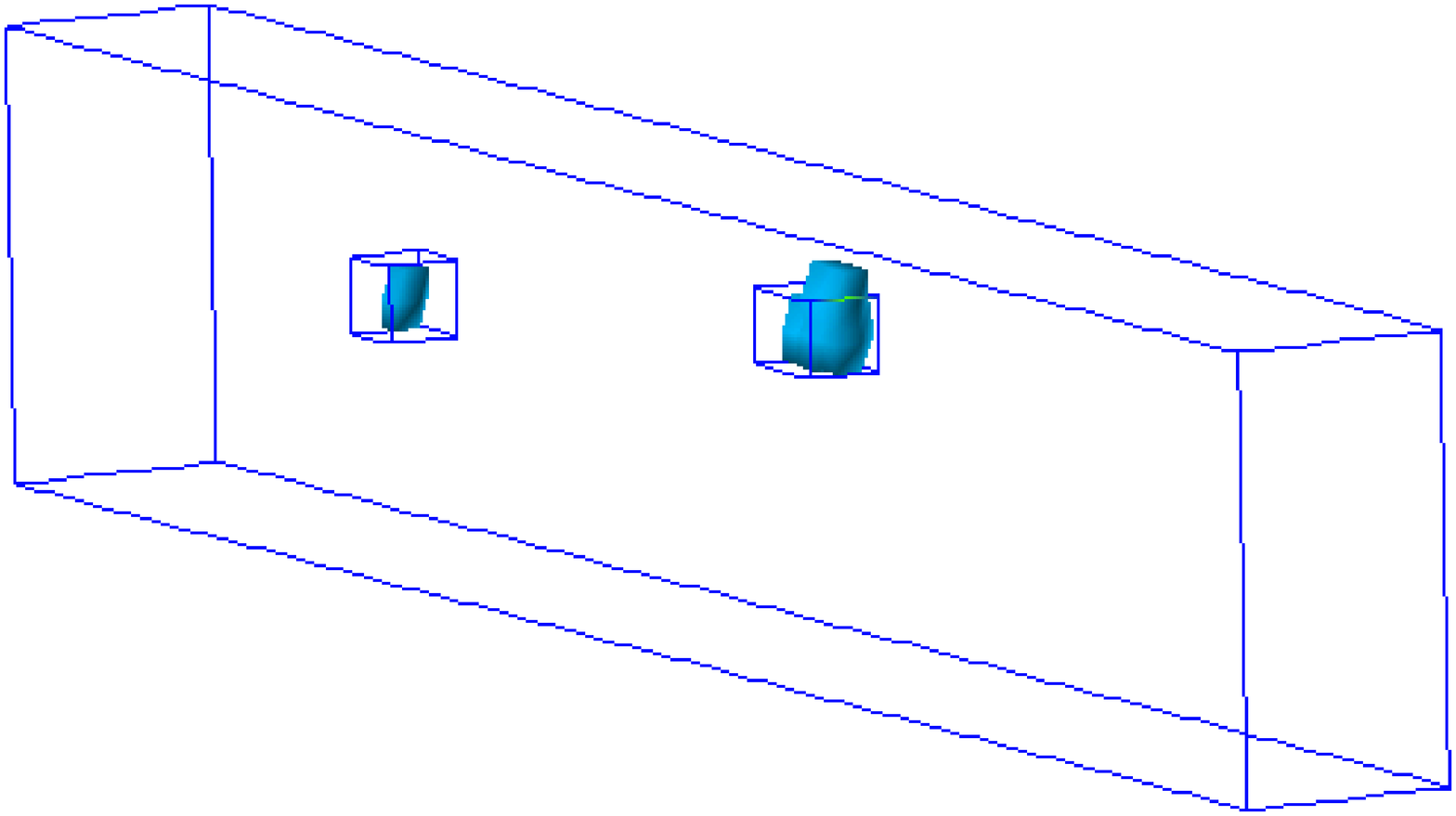}} 
 & {\includegraphics[angle=-90,width=5.5cm, clip = true]{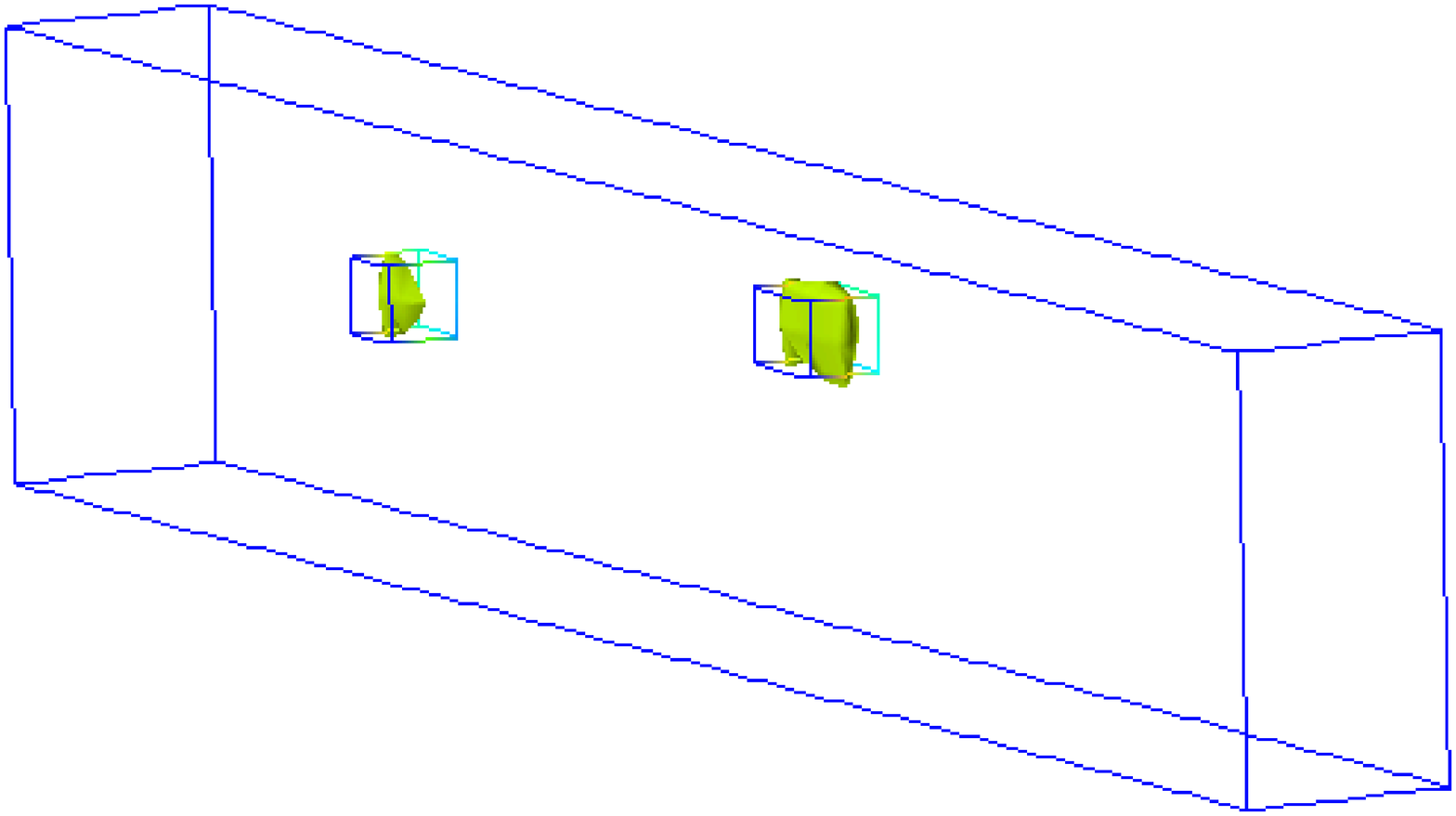}} \\
 a) $\varepsilon \approx 11.3$  &
 b) $\mu \approx 1.7$
 \end{tabular}
 \end{center}
 \caption{Computed images of  reconstructed  functions
   $\varepsilon(x)$ and $\mu(x)$ from backscattering data with 3\%
   additive noise, frequency $\omega = 30$.}
 \label{fig:Recos_omega30noise3}
 \end{figure}

\begin{figure}
 \begin{picture}(200,240)
\put(60,125) {\includegraphics[angle=-90,width=4.5cm, clip = true]{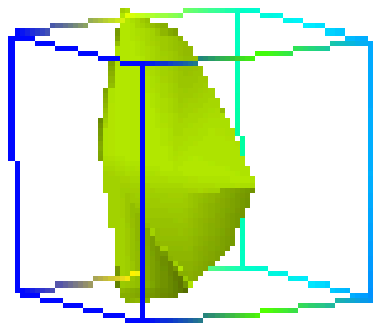}} 
\put(210,125) {\includegraphics[angle=-90,width=4.5cm, clip = true]{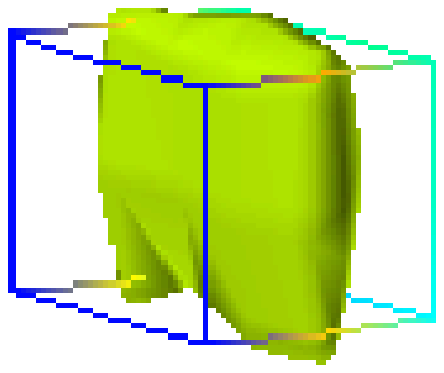}} 
\put(60,235){\includegraphics[angle=-90,width=4.5cm, clip = true]{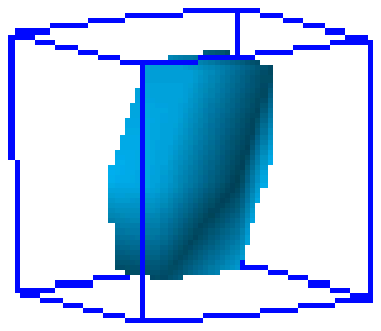}} 
\put(210,235){\includegraphics[angle=-90,width=4.5cm, clip =  true]{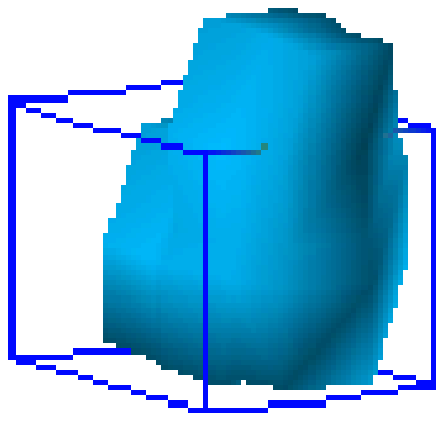}} 
 \end{picture}
\caption{Zoomed reconstructions of $\varepsilon$ (top row) and $\mu$
  (bottom row). Frequency $\omega = 30$ additive noise 3\%.  }
\label{fig:ZoomedEpsMu_omega30noise3}
\end{figure}

\begin{figure}[tbp]
 \begin{center}
 \begin{tabular}{cc}
 {\includegraphics[angle=-90,width=5.5cm, clip = true]{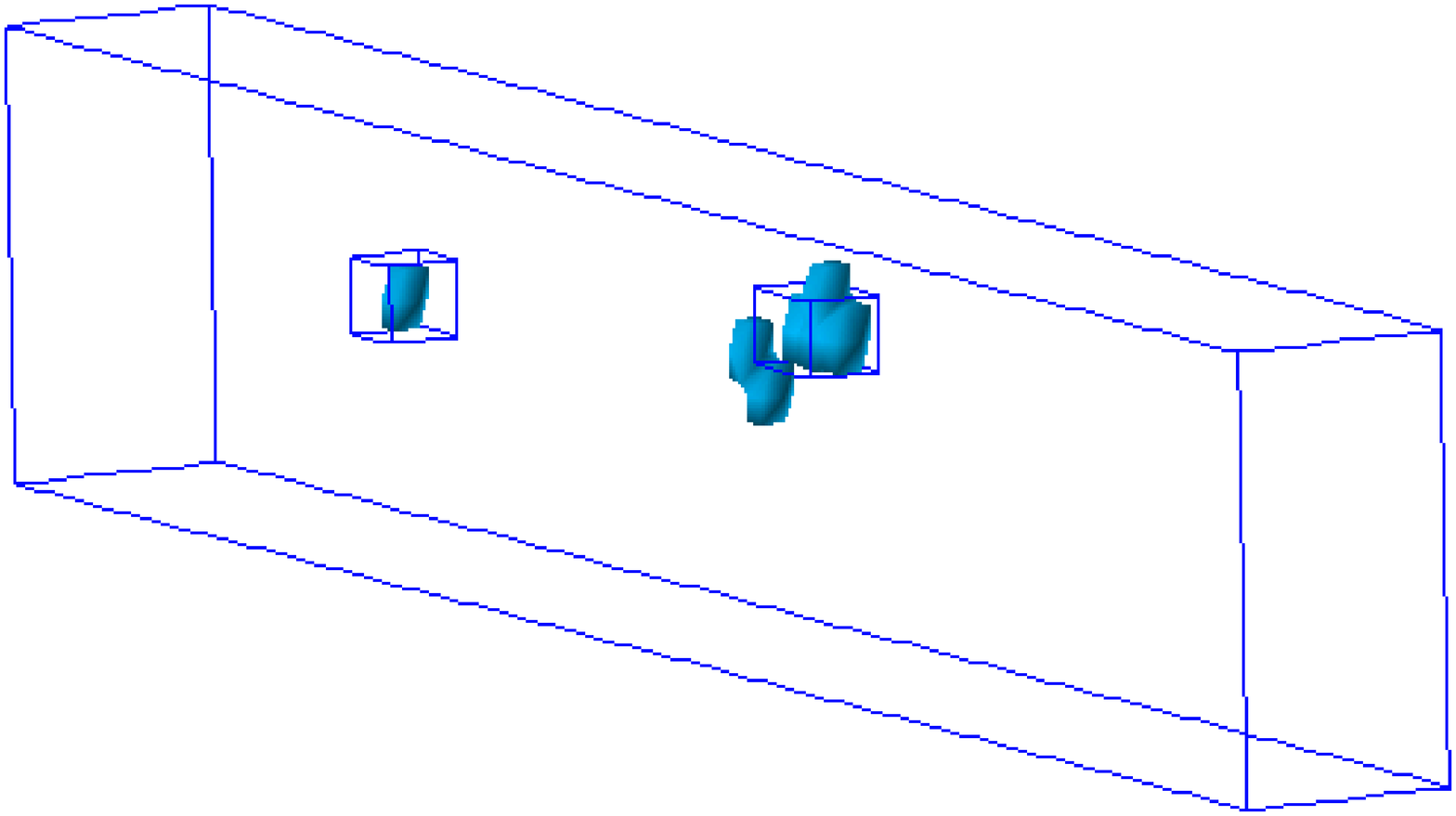}} 
 & {\includegraphics[angle=-90,width=5.5cm, clip = true]{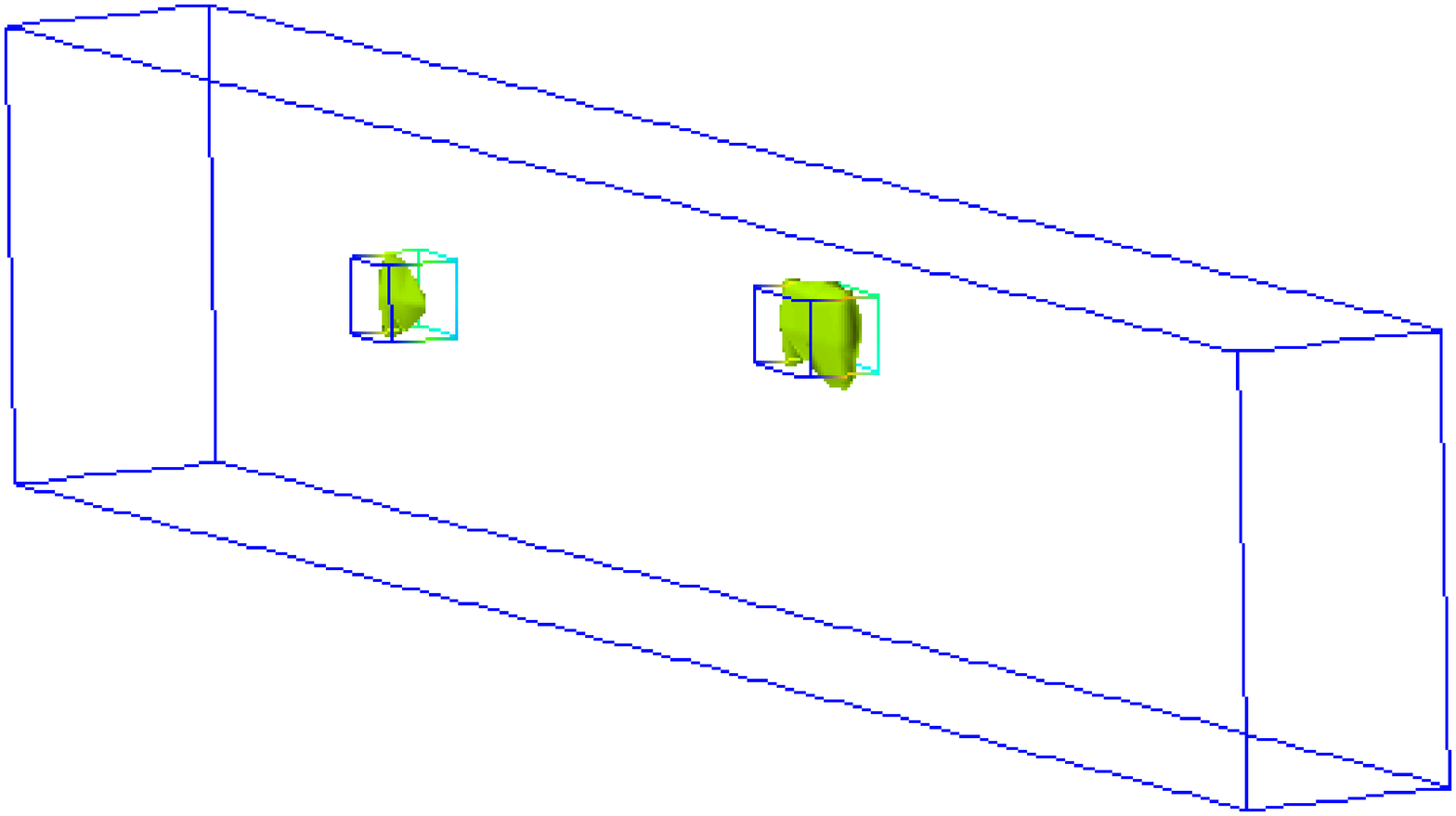}} \\
 a) $\varepsilon \approx 11.2$  &
 b) $\mu \approx 1.7$
 \end{tabular}
 \end{center}
 \caption{Computed images of  reconstructed  functions
   $\varepsilon(x)$ and $\mu(x)$ from backscattering data with 10\%
   additive noise, frequency $\omega = 30$.}
 \label{fig:Recos_omega30noise10}
 \end{figure}

 To get images of figures
 \ref{fig:Recos_omega21noise3}~-~\ref{fig:Recos_omega30noise10}, we use
 a  post-processing procedure.  Suppose that functions $\varepsilon^n(x)$
 and $\mu^l(x)$ are our reconstructions obtained by algorithm of
 section \ref{subsec:ad_alg} where $n$ and $l$ are number of iterations
 in gradient method when we have stopped to compute $\varepsilon(x)$ and
 $\mu(x)$.  Then to get images in figures
 \ref{fig:Recos_omega21noise3}~-~\ref{fig:Recos_omega30noise10}, we set
 \begin{equation*}
 \widetilde{\varepsilon }^n(x)=\left\{ 
 \begin{array}{ll}
 \varepsilon^n(x) & \text{ if }\varepsilon^n(x)>0.25 \max\limits_{\Omega }\varepsilon^n(x), \\ 
 1, & \text{ otherwise. }%
 \end{array}%
 \right.  
 \end{equation*}
 and
 \begin{equation*}
 \widetilde{\mu}^l(x)=\left\{ 
 \begin{array}{ll}
 \mu^l(x) & \text{ if }\mu^l(x)>0.87\max\limits_{\Omega }\mu^l(x), \\ 
 1, & \text{ otherwise. }%
 \end{array}%
 \right. 
 \end{equation*}

\begin{figure}
 \begin{picture}(200,240)
\put(60,125) {\includegraphics[angle=-90,width=4.5cm, clip = true]{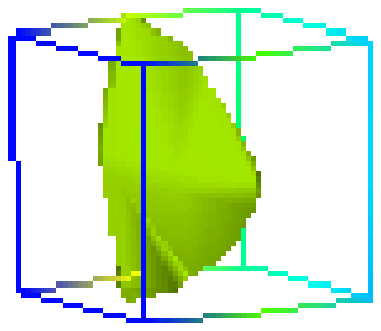}} 
\put(210,125) {\includegraphics[angle=-90,width=4.5cm, clip = true]{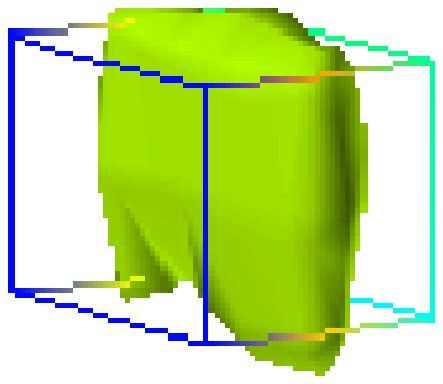}} 
\put(60,235){\includegraphics[angle=-90,width=4.5cm, clip = true]{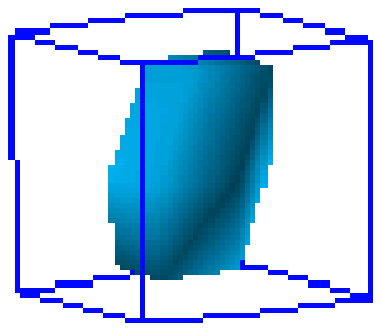}} 
\put(210,235){\includegraphics[angle=-90,width=4.5cm, clip =  true]{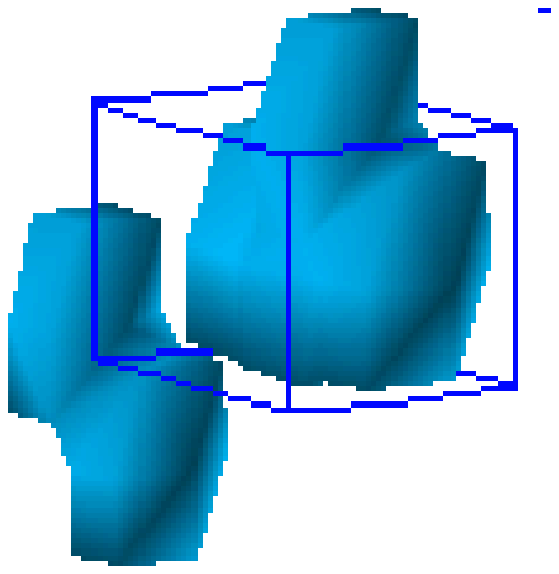}} 
 \end{picture}
\caption{Zoomed reconstructions of $\varepsilon$ (top row) and $\mu$
  (bottom row). Frequency $\omega = 30$ additive noise 10\%.  }
\label{fig:ZoomedEpsMu_omega30noise10}
\end{figure} 

\section{Discussion and Conclusion}

In this work we have used time dependent backscattering data to
simultaneously reconstruct both coefficients, $\varepsilon(x)$ and
$\mu(x)$, in the Maxwell's system as well as their locations.  In
order to do that we have used optimization approach which was similar
to the method used in \cite{BMaxwell2}.  We tested our algorithm with
two different noise levels (3\% and 10\% of additive noise) and with
two different frequencies ($\omega = 21$ and $\omega = 30$, see
(\ref{f})). The bigger noise level (10\%) seemed to produce artefact
in reconstructing $\varepsilon$ with frequency $\omega = 30$, see
figure \ref{fig:Recos_omega30noise10}. However, we are able to
reconstruct functions $\varepsilon(x)$ and $\mu(x)$ with contrasts
that are within the limits of (\ref{admpar}). An important observation
is that in our computations, we are able to obtain large contrasts for
dielectric function $\varepsilon(x)$ which allows us to conclude that we
are able to reconstruct metallic targets.
 At the same  time, the contrast
for the function $\mu(x)$ is within limits of (\ref{admpar}). We could
reconstruct size on z-direction for $\varepsilon$, however, size for
$\mu(x)$ should be still  improved.

In our future research, we are
planning to refine the obtained images through the adaptive finite element
method in order to get better shapes and sizes of the inclusions. In \cite{BTKB, BMaxwell2,  BJ} it was shown that this method is
powerful tool for the reconstruction of heterogeneous targets, their
locations and shapes accurately.

We note that our algorithm of simultaneous reconstruction of both
parameters as well as its adaptive version can also be applied for the
case when usual edge elements are used for the numerical simulation of
the forward and adjoint problems in step 1 of our algorithm, see
\cite{CWZ99, CWZ14, FJZ10} for finite element analysis and
discretization in this case.  This as well as comparison between two
different discretization techniques for the solution of our CIP can be
considered as a challenge fot the future research.

\section*{Acknowledgments}

This research was supported by the
Swedish Research Council, the Swedish Foundation for Strategic
Research (SSF) through the Gothenburg Mathematical Modelling Centre
(GMMC). Kati Niinim\"aki has been supported by the
Finnish Cultural foundation North Savo regional fund.
Part of this work was done while the second and third author were
visiting the  University of Chalmers.
 The computations were performed on
resources at Chalmers Centre for Computational Science and Engineering
(C3SE) provided by the Swedish National Infrastructure for Computing
(SNIC).

\end{document}